\documentclass[10pt,journal,compsoc]{IEEEtran}
\usepackage[utf8]{inputenc}
\usepackage[T1]{fontenc}

\usepackage{array}
\usepackage{color}

\usepackage{multirow}
\usepackage{epsfig}
\usepackage{amssymb}
\usepackage{graphicx}
\usepackage{fancyhdr}
\usepackage{multirow}
\usepackage{marginnote}
\usepackage{listings}
\usepackage{listing}

\usepackage{algorithm}
\usepackage{algorithmic}
\usepackage{xspace}

\usepackage[normalem]{ulem}
\usepackage{hyperref}

 \usepackage{dsfont}
 \usepackage{tikz}
 \usetikzlibrary{shapes.geometric}
\usetikzlibrary{arrows}

\usepackage{varwidth}

\usepackage{multirow}
 \usepackage{amsmath}

\newtheorem{definition}{Definition}[section]
\newtheorem{theorem}{Theorem}
\newtheorem{lemma}{Lemma}

\newtheorem{remark}{Remark}
\newtheorem{conjecture}{Conjecture}

\usepackage{xspace}

\newcommand{\etc}{\textit{etc.}\xspace}
\newcommand{\ie}{\textit{i.e.,}\xspace}
\newcommand{\eg}{\textit{e.g.,}\xspace}

\newcommand{\rcc}{\textsf{RC4}\xspace}
\newcommand{\ksa}{\textsf{KSA}\xspace}

\newcommand{\prga}{\textsf{PRGA}\xspace}
\newcommand{\sst}{\textsf{SST}\xspace}
\newcommand{\ctrt}{\textsl{Cyclic-to-Random Transpositions}\xspace}
\newcommand{\rtrt}{\textsl{Random-to-Random Transpositions}\xspace}
\newcommand{\riffle}{\textsl{Riffle Shuffle}\xspace}

\newcommand*{\QEDB}{\hfill\ensuremath{\square}}%

\newcommand{\E}[0]{\mathbb{E}}

\newcommand{\X}[0]{\mathbf{X}}

\newcommand{\tuple}[1]{\left \langle #1 \right \rangle}

\definecolor{ltgray}{RGB}{200,200,200}
\definecolor{dkgray}{RGB}{150,150,150}

\usepackage{float}
\newfloat{Protocol}{t}{struct}

\usepackage{framed}
\usepackage{mdwlist}
\usepackage{marginnote}
\usepackage{listings}
\usepackage{amsmath}
\usepackage{graphicx}
\usepackage{wrapfig}
\usepackage{listing}
 \usepackage[frozencache]{minted}

\ifCLASSOPTIONcompsoc
  \usepackage[nocompress]{cite}
\else
  \usepackage{cite}
\fi
\ifCLASSINFOpdf
\else
\fi
\begin{document}

%
\title{
Strong 
stationary times and its use in cryptography}
%
%
%
%
%
 
 \author{Paweł~Lorek, 
        Filip~Zagórski,
        and~Michał~Kulis
\IEEEcompsocitemizethanks{\IEEEcompsocthanksitem M. Kulis   and F. Zagórski are  with the Department of Computer Science
Faculty of Fundamental Problems of Technology, Wrocław University of Science and Technology, Poland.\protect\\
\IEEEcompsocthanksitem P. Lorek is with the Mathematical Institute, Faculty of Mathematics and Computer Science, Wrocław University, Poland. Email~Pawel.Lorek@math.uni.wroc.pl}
\thanks{All authors were supported by Polish National Science Centre contract number DEC-2013/10/E/ST1/00359.}}

\markboth{Accepted to IEEE Transactions on Dependable and Secure Computing, 
2017}%
{.}
\IEEEtitleabstractindextext{%

\begin{abstract}
This paper presents applicability of Strong Stationary Times (SST) techniques 
in the area of cryptography. The applicability is in three areas:
\begin{enumerate}
\item Propositions of a new class of cryptographic algorithms (pseudo-random 
permutation generators) which do not run for the predefined number of 
steps. Instead, these algorithms stop according to a stopping rule defined as 
SST, for which one can
obtain provable properties: 
\begin{enumerate}
\item results are perfect samples from uniform 
distribution, 
\item immunity to timing attacks (no information about the 
resulting 
permutation leaks through the information about the number of steps SST 
algorithm performed).
\end{enumerate}
\item We show how one can leverage properties of SST-based algorithms to 
construct an implementation (of a symmetric encryption scheme) which is
immune to the timing-attack by reusing 
implementations which are not secure against timing-attacks. In symmetric key 
cryptography researchers 
mainly focus on constant time (re)implementations. Our approach goes in a 
different direction and explores ideas of input masking. 
\item  Analysis of idealized (mathematical) models of existing cryptographic 
schemes -- \ie we improve a result by 
Mironov~\cite{Mironov2002}.
\end{enumerate}

\end{abstract}

}

\maketitle

\IEEEdisplaynontitleabstractindextext

%
\IEEEpeerreviewmaketitle



\IEEEraisesectionheading{\section{Introduction}\label{sec:introduction}}
 
In the traditional cryptography primitives are treated as mathematical objects,
security definitions involve 
\textsl{black-box} model: it is assumed that an adversary  has only access to 
inputs 
and outputs of a given cryptographic scheme. In reality however, a 
device, or a particular implementation, may expose lots of additional 
side-channel information,
\textsl{the running time} being one of the most important. It can be  
particularly dangerous, since one can often 
measure the time even without physical access to the device. Attacks 
exploiting this are called \textsl{timing attacks}.
\begin{algorithm}[H]
\caption{Simple modular exponentiation $c^d \bmod n$}
\label{alg:modpower}
\begin{algorithmic}[1]
\REQUIRE $c, d = d_s d_{s-1}\cdots d_1d_0; d_i \in \{0, 1\}$
\STATE $x = 1$
\FOR{i = s \TO 0} \STATE 
	{
	$x = mod(x \cdot x,  n)$
	\IF{$d_i$ == 1} \STATE {$x = mod(x \cdot c,  n)$}\ENDIF
	} 
\ENDFOR
\RETURN $x$
\end{algorithmic}
\end{algorithm}
\noindent
Straightforward implementations of RSA used the modular 
exponentiation algorithm as described in Algorithm~\ref{alg:modpower} together 
with Algorithm~\ref{alg:moddiv}.
\begin{algorithm}[H]
\caption{Modular division $mod(x, n) = x \bmod n$}
\label{alg:moddiv}
\begin{algorithmic}[1]
\REQUIRE $x, n$
\IF{$x \geq n$} \STATE 
	{ $x = x \% n$ }
\ENDIF
\RETURN $x$
\end{algorithmic}
\end{algorithm}
\noindent
An adversary, by careful 
selection of messages $c$, and just by measuring the time it takes to compute 
the result is able to completely recover $d$ -- the private key (since
the number of times line number $2$ of the Algorithm~\ref{alg:moddiv} and line 
number  $5$ of the Algorithm~\ref{alg:modpower} are called depends both on $c$ 
and $d$).
Kocher~\cite{Kocher1996} was the first one to discuss timing attacks.

How can one prevent timing attacks? There are two natural approaches:
\begin{itemize}
 \item[\textbf{A1)}]\textbf{Constant-time implementation}. Re-implementing the 
algorithm in such a way that the running time is data independent (it runs in 
the same time on all possible inputs).
 \item[\textbf{A2)}]\textbf{Masking (blinding)}. Use some extra randomness in 
such a way that time-dependent operations are performed on random inputs. 
\end{itemize}
The A1 approach is self-explaining. However, not all algorithms are easily transformable to constant time ones.
Moreover, the compiler in some cases can optimize the algorithm resulting in 
data-dependent running time  again, or a constant-time implementation on one 
platform may become again time attack susceptible on a different platform -- 
see~\cite{PorninThomas,Pereida}.

In this paper we focus our attention on masking (\ie A2) 
by performing some ``random'' permutation first, which is key-dependent.

Authors of \cite{McBits} presented fast constant-time algorithms for 
code-based public-key cryptography called \textsl{McBits}. 
McBits uses a permutation network to achieve immunity against cache-timing 
attacks. Permutation network is a network built up of wires and comparators. 
Wires carry values while comparators sort 
values between two wires (smaller value goes to the upper wires, bigger value 
goes to the lower wires).
Permutation network is a fixed construction: number of wires and comparators is 
fixed. When such a network 
has the property that it sorts all inputs, then it is called a sorting 
network.

In Section 5 of~\cite{McBits} on ``Secret permutations'' authors write (where a 
permutation is considered as the output of a permutation network)
\begin{quotation}
\textsl{``By definition a permutation network can reach every permutation, but perhaps it is much more likely
 to reach some permutations than others. Perhaps this hurts security. Perhaps not [...]''}
\end{quotation}
Above motivation leads to the following question:
\begin{quotation}
\noindent\textsl{How to obtain a key-based permutation which is 
indistinguishable from a truly random one?}
\end{quotation}
The doubt about the resulting permutation  in \cite{McBits} can come from the following:
one can have a recipe for constructing a permutation involving some number of single steps, determining 
however the necessary number of steps is not trivial. In this paper we resolve the issue, our considerations 
are Markov chain based. In particular, we consider some card shuffling schemes.
Our methods can have two applications: we can determine number of steps 
so that we are arbitrary close to uniform permutation; in so-called SST version algorithm does not run
pre-defined number of steps, but instead it stops once randomness is reached. As a side-effect we have 
  a detailed analysis of some existing stream ciphers.
Many stream ciphers (e.g., RC4,  RC4A~\cite{rc4a},   RC4+~\cite{rc4plus}, VMPC~\cite{vmpc}, Spritz~\cite{Spritz}) are 
composed of two algorithms:
\begin{enumerate}
\item \ksa  (Key Scheduling Algorithm)  uses a secret key to 
transform identity 
permutation of $n$ elements into 
some other permutation.

\item \prga (Pseudo Random Generation Algorithm) starts with a permutation 
generated by \ksa and outputs random bits from it (some function of the 
permutation and the current step number) while updating the permutation at the 
same time.
\end{enumerate}

\ksa{s} of all aforementioned algorithms (RC4, RC4A, Spritz, RC4+, VMPC) can be seen as performing some 
card shuffling, where a secret key corresponds to/replaces  randomness. 
If we consider a version of the algorithm with purely random secret key of 
infinite length, then we indeed consider a card shuffling procedure. Following 
\cite{Mironov2002}, we will refer to such a version of the algorithm 
as an \textsl{idealized version}.
In the case of \ksa used by \rcc the idealized version (mathematical model) 
of the card shuffle is called \ctrt shuffle.

The \ksa{s} of mentioned ciphers perform shuffling for some \textbf{predefined 
number of steps}.
The security of such a scheme is mainly based on analyzing  idealized version 
of the algorithm and corresponds to the ``quality of a 
shuffling''. Roughly speaking, shuffling is considered as a Markov chain on permutations, all of them converge to uniform distribution
(perfectly shuffled cards). Then we should perform as many steps as needed to be close to this uniform distribution, what is directly 
related to the so-called {\em mixing time}. This is one of the main drawbacks 
of 
\rcc: it performs \ctrt for $n$ steps, whereas the mixing time is 
of the order $n\log n$.

There is a long list of papers which point out weaknesses of the 
\rcc~algorithm. Attacks exploited both weaknesses of \prga and \ksa or the way 
\rcc~ was used 
in specific systems 
\cite{Golic1997,Fluhrer2000,Fluhrer2001a,ManSha2001,rc4royal}.
As a result, in 2015 \rcc was prohibited in TLS by IETF, Microsoft and 
Mozilla.

In the paper we use so-called \textbf{Strong Stationary Times} (SST) for 
Markov chains. The main area of application of SSTs is studying the rate of 
convergence of a chain to its stationary distribution. However, they may also be 
used for \textsl{perfect sampling} from stationary distribution of a Markov 
chain, consult \cite{Propp1996} (on   Coupling From The Past algorithm) and 
\cite{Fill1998} (on algorithms involving Strong Stationary Times and Strong 
Stationary Duality). 

This paper is an extension (and applications to timing attacks) of \cite{Zagorski_RST} presented at Mycrypt 2016 conference.

\subsection{Related work - timing attacks}

Kocher~\cite{Kocher1996} was first to observe that information about the time 
it takes a given cryptographic implementation to compute the result may 
leak information about the secret (or private) key. 
Timing attack may be treated as the most dangerous kind of side-channel 
attacks since it may be 
exploited even by a remote adversary while most of other types of side-channel 
attacks are distance bounded. 
Kocher discussed presented timing attacks on 
RSA, DSS and on symmetric key encryption schemes. 

Similar timing attack was 
implemented for CASCADE smart cards \cite{Dhem1998}.
Authors report recovering 512-bit key after 350.000 timing measurements (few 
minutes in late 90s).
Later, Brumley and Boneh~\cite{Brumley2005} showed that  
timing attacks may be performed remotely and on algorithms that were previously 
not considered vulnerable (RSA using Chinese Remainder Theorem).

The usual way of dealing with timing attacks is to attempt to re-implement a 
given scheme in such a way that it runs the same number of steps (\ie approach \textbf{A1}), no matter what 
the input is. This approach has the following drawbacks: 
\begin{enumerate}
 \item it requires to replace cryptographic hardware,
 \item it is  hard to do it right,
 \item it badly influences the performance.
\end{enumerate}

(1) The first limitation is quite obvious. So it may be worth considering the 
idea of reusing a timing attack susceptible implementation by masking (see 
Section~\ref{sec:rw-masking}) its inputs. This was 
 the original timing attack countermeasure technique proposed by Kocher for 
both RSA and DSS.

(2) An example of s2n -- Amazon's TLS implementation and its rough road to 
become immune to timing attacks~\cite{albrecht2015lucky} shows that achieving 
timing resistance is not easy. It is not an isolated case, it may be hard to 
achieve constant time in the case of of implementation~\cite{PereidaGarcia2016} 
but what may be  more worrying, even a change at the architecture 
level may not lead to success~\cite{Cachebleed}. 

(3) BearSSL is an implementation of TLS protocol which consists of many 
constant time implementations (along with a ``traditional'' ones) of the main 
cryptographic primitives. Constant-time 
implementation of CBC encryption for AES is 5.89 times slower than the 
``standard'' implementation~\cite{Pornin}. 
On a reference platform it achieves 
the following performance: AES-CBC (encryption, big): 161.40 MB/s while the 
constant time (optimized for 64-bit architecture) has 27.39 MB/s.
The slow-down is not always at the same rate, it depends on both the function 
and on the mode of operation used (for AES in CTR mode the slowdown is 1.86; 
for 3DES it is 3.10).

\subsubsection{Masking}\label{sec:rw-masking}
Let $F(\cdot)$  
be the RSA 
decryption (described as Algorithm~\ref{alg:modpower}) implementation  which 
uses a private key $[d, N]$ to decrypt an
encrypted
message $x^e \bmod N$. The following shows how to design $\mathcal{F}$ which is immune to 
timing attacks but reuses implementation of $F$.
The approach proposed by Kocher~\cite{Kocher1996} is to run a bad 
implementation $F$ on a random input (same idea as Chaum's 
blind-signatures~\cite{Chaum1983}), the code for $\mathcal{F}$ is presented in Algorithm \ref{alg:rsasecure}.

\begin{algorithm}[H]
\caption{$\mathcal{F}$ -- timing attack immune implementation of RSA decryption using 
time-attack susceptible implementation $F$}
\label{alg:rsasecure}
\begin{algorithmic}[1]
\REQUIRE $x, F, d, e, N$
\STATE select $r$ uniformly at random from $Z_N^*$
\STATE compute $y = F(d, x r^e \bmod N)$
\RETURN $s = y r^{-1} \bmod N$
\end{algorithmic}
\end{algorithm}
 
The  security of this approach comes from the fact that while the running time 
of $G$ is not constant, an attacker cannot deduce information about the 
key, because it does not know the input to the $F$. This is because $x r^e \bmod N$ is a 
random number uniformly distributed in $Z_N^*$.

The beauty of this approach lies in the fact that $F$ is actually run on a random input.
After applying $F$ one can easily get rid of it.
To be more precise: let $g_1(x, r) = x r^e \bmod N$ and let $g_2(x, r) = x r^{-1} \bmod N$ then
$G(d, x) = g_2(F(d, g_1(x, r)), r) = F(d, x)$ (both $g_1$ and $ g_2$ depend on 
the public key $[N, e]$). Taking advantage of properties of the underlying group 
lets us to remove the randomness $r$
(which remains secret to the attacker during the computations).

Similar approach can be applied to other public key schemes \eg RSA 
signatures, DSA signatures.
It is unknown how this approach -- how to design $G$ -- may be applied to 
symmetric-key cryptographic schemes without introducing any changes to the 
original function $F$.

\subsection{Our contribution}
\textbf{(1) Strong stationary time based PRNGs (Pseudo Random Number Generators).} 
Instead of running a PRNG algorithm  
for some pre-defined number of steps, we make it randomized (Las Vegas 
algorithm). As far as we are aware the idea of key-dependent running time in 
context of symmetric ciphers was not considered earlier.
To be more specific we suggest utilization of so-called \textbf{Strong 
Stationary Times} (SST) for Markov chains.
We use SST to 
obtain samples from uniform distribution on all permutations (we actually 
perform perfect sampling). 
Few things needs pointing out:
  \begin{enumerate}
  \item In McBits \cite{McBits} construction, authors perform \textsl{similar} masking
  by performing some permutation. This is achieved via a permutation network.
  However, they: a) run it for a pre-defined number of steps (they are  unsure if the resulting permutation
  is random); b) they work hard to make the implementation constant-time.
  \par 
  Utilization of SST allows to obtain a \textsl{perfect} sample from uniform distribution (once SST rule is found, 
  we have it ``for free''). This is a general concept which may be applied to 
many algorithms.
  \item Coupling methods are most commonly used 
tool for studying the rate of 
convergence to stationarity for Markov chains.
   They allow to bound the so-called \textsl{total variation distance} between the distribution of given chain at time instant $k$ and 
   its stationary distribution. However, the (traditional) security definitions 
require something ``stronger''. It turns out that bounding
   \textsl{separation distance} is what one actually needs. It fits perfectly 
into the notion of Strong Stationary Times we are using (see 
Section~\ref{sec:sst_security}).
   
  \item By construction, the running time of our model is key dependent. 
In extreme cases (very unlikely) it may leak some information about the key, but
  it \textbf{does not leak any} information about the resulting permutation to 
an adversary. We elaborate more on this in  Section \ref{app:timing}.

  \end{enumerate}

\textbf{(2) Preventing timing attacks.}
We present a technique which allows to obtain a timing attack immune 
implementation of a symmetric key encryption scheme by re-using a timing 
attack vulnerable scheme. The technique follows ideas of Kocher's masking of 
the RSA algorithm and are similar in nature to the utilization of a permutation network in 
McBits.

The idea is that an input to a timing-attack vulnerable function is first a 
subject to a pseudo-random permutation. The permutation is generated (from a 
secret key) via SST-based pseudo-random permutation generation algorithm.

The main properties of such an approach are following:
the resulting permutation is uniformly random provided the long enough key is 
applied; the running time of an SST-based algorithm does not reveal any 
information about the resulting permutation. 

Based on our straightforward implementation,
we show that the performance of our approach is  similar to constant-time 
implementations
available in BearSSL~\cite{Pornin}, see Section~\ref{sec:conclusions}, Figure 
\ref{tab:efficiency} (but does not require that much effort and caution 
to implement it right).

\textbf{(3) Analysis of idealized models of existing cryptographic schemes.}
 
The approach, in principle, allows us to analyze the quality of existing PRNGs, 
closing the gap between theoretical models and practice. The ``only'' thing 
which must be done for a given scheme is to 
find an SST for the corresponding shuffle and to analyze its properties. We 
present in-depth analysis of  \rcc's \ksa algorithm.
 As a result of Mironov's~\cite{Mironov2002} 
work, one knows that the
idealized version of \rcc's \ksa would need keys of length $\approx 23~037$ in 
order to meet the mathematical model. 
 Similarly as in \cite{Mironov2002}, we propose  SST which is valid for \ctrt and \rtrt.
 For the latter one we calculate the  mixing time  which is ``faster''.
It is known that \rtrt card shuffling needs ${1\over 2}n\log n$ steps to mix. 
It is worth mentioning that this shuffling scheme  and \ctrt are  similar in the \textsl{spirit}, it 
does not transfer automatically  that the latter one also needs ${1\over 2}n\log n$  steps to mix. 
Mironov  \cite{Mironov2002} states that the  expected time till reaching uniform 
distribution is upper bounded by $2n H_n - n$, we show that the 
expected running time for this SST applied to \rtrt is equal to 
$$ E[T]  =   nH_n + n + O(H_n)$$
($H_n$ is the $n$-th harmonic number)
and empirically check  that the result is 
similar for \ctrt.
This directly translates into the required steps that should be performed by 
\rcc.

In fact one may use a better 
shuffling than the one that is used in \rcc, \eg
\textit{time-reversed riffle shuffle} which requires (on average) $4 096$ bits 
-- not 
much more than $2 048$ bits which are allowed for \rcc (see 
Section~\ref{sec:nasz_ksa}).


\section{Preliminary}
Throughout  the paper, let $\mathcal{S}_n$ denote a set of all 
permutations 
of a set $\{1,\ldots,n\}=:[n]$.

\subsection{Markov chains and rate of convergence}\label{sec:mcRate}
Consider an ergodic Markov chain $\X=\{X_k, k\geq 0\}$ on a finite state space 
$\E=\{0,\ldots,M-1\}$ with a stationary distribution $\psi$. Let $\mathcal{L}(X_k)$ denote the distribution of the chain at time instant $k$. 
By the  \textsl{rate of convergence} we understand the knowledge on how fast the distribution of the chain converges to
its stationary distribution. We have to measure it according to some distance $dist$. 
Define
$$\tau^{dist}_{mix}(\varepsilon)=\inf\{k: dist(\mathcal{L}(X_k),\psi)\leq \varepsilon\},$$
which is called \textsl{mixing time} (w.r.t. given distance $dist$).
In our case the state space is the  set of all permutations of $[n]$, \ie  $\E:=\mathcal{S}_n$.
The stationary distribution is the uniform distribution over $\E$, \ie  $\psi(\sigma)={1\over n!}$ for all $\sigma\in\E$.
In most applications the mixing time is defined w.r.t. total variation distance:
$$d_{TV}(\mathcal{L}(X_k),\psi)={1\over 2} \sum_{\sigma\in \mathcal{S}_n} \left| Pr(X_k=\sigma)-{1\over n!}\right|.$$
The \textbf{separation distance} is defined by:
$$\begin{array}{lll}
sep(\mathcal{L}(X_k),\psi):=\displaystyle\max_{\sigma\in\E} \left(1-n!\cdot Pr(X_k=\sigma)\right). \\
\end{array}
$$
It is relatively easy to check that $d_{TV}(\mathcal{L}(X_k),\psi)\leq sep(\mathcal{L}(X_k),\psi)$.

\textbf{Strong Stationary Times.} The definition of separation distance fits 
perfectly into the
notion of Strong Stationary Times (SST) for Markov chains.
This  is a probabilistic tool for studying the rate of convergence of Markov 
chains allowing also \textsl{perfect sampling}.
We can think of stopping time as of a running time  of an algorithm which 
observes a
Markov chain $\X$
and which stops according to some stopping rule (depending only on the past).
\begin{definition}
 Random variable $T$ is a randomized stopping time if it is a running time of 
the \textsc{Randomized Stopping Time} algorithm.
 
\begin{algorithm}[H]
\caption{\textsc{Randomized Stopping Time}}\label{euclid}
 \begin{algorithmic}[1]\label{alg:rand_stop_time}
 \STATE $k:=0$
 \STATE $coin:=$\texttt{Tail}
 \WHILE{ $coin==$\texttt{Tail}}
  \STATE At time $k, (X_0,\ldots,X_k)$ was observed 
  \STATE Calculate $f_k(X),$ where $f_k: (X_0,\ldots,X_k) \to [0,1]$
  \STATE Let $p=f_k(X_0,\ldots,X_k)$. Flip the coin resulting in \texttt{Head} 
with probability 
  $p$ and in \texttt{Tail} with probability $1-p$. Save result as $coin$.
  \STATE $k:=k+1$
\ENDWHILE
\end{algorithmic}
\end{algorithm}
\end{definition} 

\begin{definition}\label{def:SST}
 Random variable $T$ is a \textbf{Strong Stationary Time (SST)} if it is a randomized stopping time  for chain $\X$ such that:
 $$\forall(i\in \E)\ Pr(X_k=i | T=k) = \psi(i).$$ 
\end{definition}

Having SST $T$   lets us bound the following 
(see \cite{Aldous1987})
\begin{equation}\label{eq:TV_sep_SST}
sep(\mathcal{L}(X_k),\psi)\leq Pr(T>k).
\end{equation}

\subsection{Distinguishers and security definition}
In this section we define a \textsf{Permutation 
distinguisher} -- we 
consider its advantage in a  traditional cryptographic security definition.
That is, a \textsf{Permutation distinguisher}  is given a permutation $\pi$ and 
needs to decide whether $\pi$ is a result of a given (pseudo-random) algorithm 
or if $\pi$ is a permutation 
selected uniformly at random from the set 
of all permutations of  a given size. The upper bound on \textsf{Permutation 
distinguisher}'s advantage is
an upper bound on {\underline{any}} possible distinguisher -- even those which
are not bounded computationally. Additionally, in the 
Section~\ref{sec:signDist} we consider  \textsf{Sign distinguisher}.

\textbf{Permutation distinguisher} 
The distinguishability game for the adversary is as follows 
(Algorithm~\ref{def:shuffle}):
\begin{definition}{The permutation indistinguishability \textsf{Shuffle$_{\textsf{S},\mathcal{A}}(n, r)$} experiment:} 
\begin{algorithm}[H]
\caption{\textsf{Shuffle$_{\textsf{S},\mathcal{A}}(n, 
r)$}}\label{def:shuffle}
Let $S$ be a shuffling algorithm which  in each round requires $m$ bits,
and let $\mathcal{A}$ be an adversary.
\begin{enumerate}
    \item \textsf{S} is initialized with:
      \begin{enumerate}
      \item a key generated uniformly at random $\mathcal{K} \sim 
\mathcal{U}(\{0, 1\}^{r m})$,
      \item $S_0 = \pi_{0}$ (identity permutation)
      \end{enumerate}
  \item \textsf{S} is run for $r$ rounds: $S_r := \textsf{S}(\mathcal{K})$ and 
produces a permutation $\pi_r$. 
 \item We set: 
  \begin{itemize}
    \item  $c_0 := \pi_{rand}$ a random permutation from uniform 
      distribution is chosen,
    \item $c_{1} := \pi_r.$
    \end{itemize}
    \item A challenge bit $b \in \{0, 1\}$ is chosen at random, permutation 
$c_b$ is sent to the Adversary.
  \item Adversary replies with $b'$.
  \item The output of the experiment is defined to be $1$ if $b' = b$, and $0$ 
otherwise.
\end{enumerate}
\end{algorithm}

\end{definition}
In the case when the adversary wins the game (if $b = b'$) we say that 
$\mathcal{A}$ succeeded. The adversary wins the game if 
it can distinguish the random permutation from the permutation being a result 
of the 
PRPG algorithm.

\begin{definition}\label{def:indist}
 A shuffling algorithm \textsf{S} generates indistinguishable 
permutations 
if for all adversaries $\mathcal{A}$ 
there exists a negligible function $\textsf{negl}$ such that 

\[
 Pr\left[\textsf{Shuffle$_{\textsf{S},\mathcal{A}}(n,r)$} = 
1\right] \leq 
\frac{1}{2} + \textsf{negl}(n).
\]
\end{definition}
The above translates into:
\begin{definition}\label{def:advantage}
 A shuffling algorithm \textsf{S} generates indistinguishable 
permutations 
if for any adversary $\mathcal{A}$ 
there exists a negligible function $\textsf{negl}$ such that:

\begin{eqnarray*}
 \left|
\underset{\scriptscriptstyle K \leftarrow \{0, 1\}^{keyLen}}{Pr}
 \left[\mathcal{A}( \textsf{S}(K)) = 1 \right] -  
{\underset{\scriptscriptstyle R \leftarrow 
\mathcal{U}(\mathcal{S}_n)}{Pr}}[\mathcal{A}( 
R) = 1]   
\right|  \\ 
\leq \textsf{negl}(n).
\end{eqnarray*}

\end{definition}

\subsection{Strong stationary times and security guarantees}\label{sec:sst_security}

Coupling method is a commonly used tool for bounding the rate of convergence 
of Markov chains. 
Roughly speaking, a coupling of a Markov chain $\X$ with a transition matrix 
$\mathbf{P}$ is a 
bivariate chain $(\X',\X'')$ such that 
 marginally  $\X'$ and $\X''$ are Markov chains with the transition matrix $\mathbf{P}$ and once the chains meet they stay together
(in some definitions this condition  can be relaxed). Let then 
$T_c=\inf_k\{X'_k=X''_k\}$ \ie the first time chains meet, called 
\textsl{coupling time}. The  \textsl{coupling inequality} states that 
$d_{TV}(\mathcal{L}(X_k),\psi)\leq Pr(T_c>k)$. 

On the other hand, separation distance is an upper bound on total variation 
distance, \ie $d_{TV}(\mathcal{L}(X_k),\psi)\leq sep(\mathcal{L}(X_k),\psi)$.
At first glance it seems that it is better to directly bound $d_{TV}$, since we can have $d_{TV}$ very small, whereas $sep$ is (still) large.
However, knowing that $sep$ is small gives us much more than just knowing that 
$d_{TV}$ is small, what turns out to be \textbf{crucial}
for proving security guarantees (\ie Definition~\ref{def:advantage}). In our 
case ($\E=\mathcal{S}_n$ and $\psi$ is a uniform distribution on $\E$)
having $d_{TV}$ small, \ie $d_{TV}(\mathcal{L}(X_k),\psi)={1\over 2}\sum_{\sigma\in\mathcal{S}_n}\left|Pr(X_k=\sigma)-{1\over n!}\right|\leq \varepsilon$ 
\textbf{does not} imply that  $|Pr(X_k=\sigma)-{1\over n!}|$ is uniformly small
(\ie  of order ${1\over n!}$). Knowing however  that $sep(\mathcal{L}(X_k),\psi)\leq \varepsilon$ 
implies \begin{equation}\label{eq:sep_unif}
 \forall(\sigma\in\E) \ \left|Pr(X_k=\sigma)-{1\over n!}\right|\leq {\varepsilon\over n!}.
\end{equation}
Above inequality is what we need in our security definitions and shows that the notion of separation distance is an adequate measure 
of mixing time for our applications.

It is worth noting that  $d_{TV}(\mathcal{L}(X_k),\mathcal{U}(\E))\leq 
\varepsilon$ implies (see Theorem~7 in 
\cite{Aldous1986})
 that $sep(\mathcal{L}(X_{2k}),\mathcal{U}(\E))\leq \varepsilon$. This  means 
that proof of security which bounds directly total variation distance by 
$\varepsilon$ would require twice as many 
bits of randomness compared to the result which guarantees $\varepsilon$ bound 
on the separation distance.

\section{Shuffling cards}\label{sec:shuffling}
In this Section we consider four chains $\X^{(T2R)}, \X^{(RTRT)}, \X^{(CTRT)}, \X^{(RS)}$
corresponding to known shuffling schemes: \textsl{Top-To-Random}, \textsl{Random-To-Random Transpositions}, \textsl{Cyclic-To-Random Transpositions} 
and \textsl{Riffle Shuffle} (we actually consider its time reversal). The state space of all chains is $\E=\mathcal{S}_n$.
The generic state is denoted by $S\in\E$, where $S[j]$ denotes the position of element $j$.
Rough description of one step of each shuffling (in each shuffling ``random'' corresponds to the uniform distribution) is following:
\begin{itemize}
 \item $\X^{(T2R)}$: choose random position $j$ and move the card $S[1]$ (i.e., which is currently on top)  to this position.
 \item $\X^{(RTRT)}$: choose two positions $i$ and $j$ at random and swap cards $S[i]$ and $S[j]$.
 \item $\X^{(CTRT)}$: at step $t$ choose random position $j$ and swap cards $S[(t \bmod n)+1]$ and $S[j]$.
 \item  $\X^{(RS)}$: assign each card random bit. Put all cards with bit assigned 0 to the top keeping relative positions of these cards.
\end{itemize}
The stationary distribution of all the chains is uniform, i.e. $\psi(\sigma)={1\over n!}$ for all $\sigma\in\E$.
Aforementioned Strong Stationary Time is one of the methods to study rate of convergence.
Probably the simplest (in description) is an SST for $\X^{(T2R)}$: 
\begin{quote}
Wait until the card which is initially at the bottom
reaches the top, then make one step more.
\end{quote}
Denote this time by $T^{(T2R)}$. In this classical example $T^{(T2R)}$ corresponds to the classical coupon collector problem.
We know that  $Pr(T^{(T2R)}>n\log n+cn)\leq e^{-c}$, thus via (\ref{eq:TV_sep_SST}) we have 
$$d_{TV}(\mathcal{L}(X^{(T2R)}_k),\psi)\leq sep(\mathcal{L}(X^{(T2R)}_k),\psi)\leq e^{-c}$$
for $k=n\log n+cn$. This bound is actually tight, in a sense that there is a so-called cutoff at $n\log n$, see Aldous and Diaconis \cite{Aldous1986}.
\par 
The chain $\X^{(RTRT)}$ mixes at ${1\over 2} n\log n$, the result was first proven by Diaconis and Shahshahani \cite{Diaconis1981a}
with analytical methods, Matthews \cite{Matthews1988} provided construction of SST for this chain. 
\par 
There is a simple construction of SST for $\X^{(RS)}$: initially consider all pairs of cards as \textsl{unmarked}.
Then, if at some step cards in some pair were assigned different bits \textsl{mark} this pair.  Stop at the 
first time when all pairs are marked. This SST yields in $2n\log_2 n$ mixing time (the actual answer is that $(3/2)n \log_2 n$ steps are enough).
\par 
The case $\X^{(CTRT)}$ is a little bit more complicated. Having SST $T$ does not automatically imply that 
we are able to bound $Pr(T>k)$, or even calculate $ET$. This is due to the cyclic structure of the shuffle. Note that \textsl{``in spirit''} 
$\X^{(CTRT)}$ is similar to $\X^{(RTRT)}$, which - as mentioned - mixes at ${1\over 2}n\log n$ (and even more - there is a so-called cutoff at
this value, meaning that ``most of the action happens'' around ${1\over 2} n\log n$). 
Mironov \cite{Mironov2002} (in context of analysis of RC4) showed that it takes $O(n\log n)$ for this shuffle to mix.
Unfortunately Mironov's ``estimate of the rate of growth of the strong stationary
time $T$ is quite loose'' and results ``are a far cry both from the provable 
upper and lower bounds on the convergence rate''. He:
 \begin{itemize}
  \item  proved an upper bound $O(n\log n)$. More precisely Mironov showed that 
there exists some positive constant 
  $c$ such that $P[T>c n\log n]\to0$ when $n\to\infty.$ Author experimentally 
checked that $P[T > 2 n \lg n]  < 1/n$ for $n=256$ which corresponds 
to $P[T > 4096] < 1/256$. 
\item experimentally showed that $E[T] \approx 11.16n \approx 1.4 n \lg n 
\approx 2857$ (for $n = 256$) -- which 
translates into: on average one needs to drop $\approx 2601$ initial bytes.
 \end{itemize}
Later Mosel, Peres and Sinclair~\cite{Mossel2004} proved a matching lower 
bound establishing mixing time to be of order $\Theta(n\log n)$. However, the 
constant was  not determined.
\medskip\par 
One can apply the following approach: find $T$ which is a proper SST for \textsl{both} $\X^{(CTRT)}$ and $\X^{(RTRT)}$, calculate 
its expectation for $\X^{(RTRT)}$ and show experimentally that the expectation is similar for $\X^{(CTRT)}$. 
This is the main part of this Section resulting in a conjecture that the mixing time of $\X^{(CTRT)}$ is $n\log n$. \par 

\subsection{Mironov's SST  for \ctrt}
Let us recall Mironov's \cite{Mironov2002} construction of SST:
\begin{quotation}
 ``At the beginning all cards numbered $0,\ldots,n-2$ are \textsl{unchecked}, the $(n-1)^{th}$ card is \textsl{checked}.
 Whenever the shuffling algorithm exchanges two cards, $S[i]$ and $S[j]$, one of the two rules may apply before the swap takes place:
 \begin{itemize}
  \item[a.] If $S[i]$ is unchecked and $i=j$, check $S[i]$.
  \item[b.] If $S[i]$ is unchecked and $S[j]$ is checked, check $S[i]$.
 \end{itemize}
  The event $T$ happens when all cards become checked. ''
\end{quotation}
Then the author proves that this is a SST for \ctrt and shows that there exists constant $c$ 
(can be chosen less than 30) such that $Pr[T>c n\log n]\to 0$ when $n\to \infty$.
Note that this marking scheme is also valid for \rtrt shuffling, for which we have:
\begin{lemma}
The expected running time of \rtrt shuffling with Mironov stopping rule is: $$ET 
= 2nH_n -n + O(H_n).$$
\end{lemma}
\textsl{Proof.}
We start with one card checked. When $k$ cards are checked, then probability of 
checking another one is equal to $p_k = \frac{(n-k)(k+1)}{n^2}$. Thus, the time 
to check all the cards is distributed as a sum of geometric random variables 
and its expectation is equal to:
$$\sum_{k=1}^{n-1} \frac{1}{p_k} = 2 \frac{n^2}{n+1} H_n - n = 2nH_n - n + 
O(H_n).$$
 \QEDB
\subsection{Better SST  for \ctrt}\label{sec:betterSST}
We suggest another ``faster'' SST which is valid for both \ctrt and \rtrt. We 
will calculate its expectation and variance for \rtrt and check 
experimentally 
that it is similar if the 
stopping rule is applied to \ctrt. 
\noindent
The SST  is given in    \texttt{StoppingRuleKLZ} algorithm.
 \begin{algorithm}[H]  
\begin{algorithmic}[htb!]  
  \caption{\texttt{StoppingRuleKLZ}}
      \STATE \textbf{Input} set of already marked cards $\mathcal{M}\subseteq\{1,\ldots,n\}$, round $r$, $Bits$
      \STATE \textbf{Output} \{\texttt{YES,NO}\}     
      \STATE \
      \STATE $j=$\texttt{n-value}($Bits$)
      \IF {there are less than $\lceil(n-1)/2\rceil$ marked cards}
		\IF {both $\pi[r]$ and $\pi[j]$ are unmarked}
		  \STATE mark card $\pi[r]$
		  \ENDIF		  
         \ELSE
         \IF {($\pi[r]$ is unmarked and $\pi[j]$ is marked) OR ($\pi[r]$ is unmarked and $r=j$)}
	      \STATE mark card $\pi[r]$
	      \ENDIF	
\ENDIF
\STATE \
      \IF { all cards are marked}
        \STATE \texttt{STOP} 
        \ELSE 
        \STATE \texttt{CONTINUE}
        \ENDIF
\end{algorithmic}
\end{algorithm}        
 
 \begin{lemma}
 The stopping rule \texttt{StoppingRuleKLZ} is SST for $\X^{(CTRT)}$.
 \end{lemma}
\textsl{Proof.}
First phase of the procedure (\ie the case when there are less than $\lceil(n-1)/2\rceil$ cards marked) 
is a construction a random permutation of marked cards by  
placing unmarked cards on randomly chosen unoccupied positions -- this is actually the 
first part of Matthews's \cite{Matthews1988} marking  scheme.
Second phase is simply a  Broder's construction. 
Theorem~9. in \cite{Mironov2002} shows that this is a valid SST for \ctrt.
Both phases combined produce a random permutation of all cards.
 \QEDB
 
 \begin{remark}
  One important remark should be pointed. Full Matthews's marking \cite{Matthews1988} scheme is an SST which is ``faster'' than ours.
  However, although it is valid  for \rtrt, it  is not valid  for \ctrt.
 \end{remark}
Calculating $ET$ or $VarT$ seems to be a challenging task.
But note that marking scheme \texttt{StoppingRuleKLZ} also yields a valid SST  for \rtrt.
In next Lemma we calculate $ET$ and $VarT$ for this shuffle, later we experimentally show that 
$ET$ is very similar for both marking schemes.

 \begin{lemma}\label{lem:ksa_ET}
  Let $T$ be the SST   \texttt{StoppingRuleKLZ} applied to \rtrt. Then we have
  \begin{equation}\label{eq:ksa_ET}
  \begin{array}{lll}
   E[T] & = &  nH_n + n + O(H_n), \\[8pt]
  Var[T] & \sim & \frac{\pi^2}{4}n^2,    
  \end{array}  
  \end{equation}
where $f(k) \sim g(k)$ means that $\lim_{k \to \infty} \frac{f(k)}{g(k)} = 1$.
   \end{lemma}

\textsl{Proof}.
 Define  $T_k$ to be the first time when $k$ cards are marked (thus $T\equiv T_n$). Let $d=\lceil(n-1)/2\rceil$. Then $T_d$ is the running time of the first phase and $\left(T_n - T_d\right)$ is the running time of 
 the second phase. Denote $Y_k := T_{k+1}-T_k$.

Assume that there are $k<d$ marked cards at a certain step. Then the new card will be marked in next step if we choose two unmarked cards what happens with probability: 
$$p_a(k)= {(n-k)^2\over n^2}.$$
Thus $Y_k$ is a geometric random variable with parameter $p_a(k)$ and
\smallskip\par \noindent
$E[T_d] =$
$$
\begin{array}{lllllllll}
  =  \displaystyle \sum_{k=0}^{d-1} E[Y_k] = \sum_{k=0}^{d-1} {1\over p_a(k)}=\sum_{k=0}^{d-1}  {n^2\over (n-k)^2} = n^2 \sum_{k=n-d+1}^n{1\over k^2} \\[14pt] 
    = \displaystyle   n^2\left(H_n^{(2)}-H_{n-d}^{(2)}\right) = n^2 \left( \frac{1}{n} + O\left(\frac{1}{n^2}\right) \right) = n + O(1).
\end{array}
$$
%
Now assume that there are $k\geq d$ cards marked at a certain step. Then, the new card will be marked in next step with probability:
$$p_b(k)= {(n-k)(k+1)\over n^2}$$
and $Y_k$ is a geometric random variable with parameter $p_a(k)$. Thus:
\smallskip\par \noindent
$E[T_n-T_d]$
$$
\begin{array}{lllllllll}
  =    \sum_{k=0}^{d-1} E[Y_k] = \sum_{k=d}^{n-1} {1\over p_b(k)}=\sum_{k=d}^{n-1} {n^2\over (n-k)(k+1)} \\[8pt]
  =    {n^2 \over n+1} \sum_{k=d}^{n-1} \left({1\over n-k} + {1\over k+1}\right)={n^2 \over n+1} \left( {\sum_{k=1}^{n-d} \frac{1}{k}} + {\sum_{k=d+1}^{n} \frac{1}{k}} \right) \\[8pt] 
  =      {n^2\over n+1} \left( H_{n-d} + H_n - H_d \right) = {n^2\over n+1} H_n + {n^2\over n+1} \left( H_{n-d} - H_d \right) \\[8pt]
  =  nH_n - \frac{n}{n+1}H_n + {n^2\over n+1} \left( H_{n-d} - H_d \right) \\[8pt]
  = nH_n + O(H_n) + O(1) = nH_n + O(H_n).
\end{array}
$$

\noindent 
For variance we have: 
\medskip\par\noindent 
$$
\begin{array}{lllllllll}
Var[T_d]=\\[8pt]
   \displaystyle \sum_{k=0}^{d-1} \frac{1-p_a(k)}{\left(p_a(k)\right)^2} = \sum_{k=0}^{d-1} \frac{1-\frac{(n-k)^2}{n^2}}{\left(\frac{(n-k)^2}{n^2}\right)^2} \approx \int_{0}^{\frac{n}{2}} \frac{1-\frac{(n-x)^2}{n^2}}{\left(\frac{(n-x)^2}{n^2}\right)^2} dx  = \frac{4}{3} n. \\[18pt]
  Var[T_n-T_d]=\\[8pt]
  \quad \displaystyle  \sum_{k=d}^{n-1} Var[Y_k] = \sum_{k=d}^{n-1} \frac{1-p_b(k)}{\left(p_b(k)\right)^2} = \sum_{k=d}^{n-1} \frac{1 - \frac{(n-k)(k+1)}{n^2}}{\frac{(n-k)^2 (k+1)^2}{n^4}} \\[8pt]
 =  \displaystyle n^2  \sum_{k=d}^{n-1} \frac{n(n-1)+k(1-n)+k^2}{(n-k)^2 (k+1)^2}  \\[8pt]
 \approx \displaystyle n^2 \cdot \frac{1}{2} \sum_{k=0}^{n-1} \frac{n(n-1)+k(1-n)+k^2}{(n-k)^2 (k+1)^2} \\[8pt] 
  \approx  \displaystyle \frac{n^2}{2} \left[ \left(\frac{2}{n} - \frac{4}{n^2} \right) H_n + 3 H_n^{(2)} \right] \sim \frac{\pi^2}{4}n^2. \\
\end{array}
$$
\noindent Finally,
 $$Var[T_n] = Var[T_d]+Var[T_n-T_d] \sim \frac{\pi^2}{4}n^2. $$
\QEDB

\subsection{Experimental results}\label{sec:expResults}

The expectation of the following SSTs applied to \rtrt shuffling scheme is known:
\begin{itemize}
  \item Mironov's SST (used in \cite{Mironov2002}) it is $2n H_n-n$
  \item \texttt{StoppingRuleKLZ} SST  it is $n(H_n+1)$
\end{itemize}
As shown, both marking schemes are valid  for both, \ctrt and \rtrt shufflings.
For both stopping rules applied to \ctrt no precise results on expected running 
times are known. Instead we estimated them via simulations, simply running 
$10.000$ of them. The results are given in Fig. \ref{fig:sims}. 

\begin{figure}[H]
\begin{tabular}{lll}
\begin{tabular}[t]{|l|l|c|c|c|c|c|}\hline
&  &  \texttt{\small  StoppingRuleKLZ}  &  \small{Mironov's} \texttt{ST}   \\\hline
\multirow{3}{*}{\small $ET$} & $n=256$ & 1811   & 2854\\
& $n=512$ & 3994  & 6442 \\
 &    $n=1024$ & 8705  & 14324 \\ \hline 

 \multirow{3}{*}{\small $VarT$} & $n=256$ & 111341  & 156814 \\
& $n=512$ & 438576 &   597783  \\
 &    $n=1024$ & 1759162   & 2442503  \\ \hline 

\end{tabular}
\\
\begin{tabular}[t]{|c|c|c|c|c|c|c|}\hline
&  \   $n(H_n+1)$  \ &  \ $2n H_n - n$  \  \\\hline
$n=256$  & 1823.83 & 2879.66 \\
$n=512$  & 4002.05 & 6468.11 \\
$n=1024$  & 8713.39 & 14354.79 \\ \hline 
\end{tabular}
\end{tabular}\caption{Simulations' results for Mironov's and \texttt{StoppingRuleKLZ} stopping rules.}\label{fig:sims}
\end{figure}
The conclusion is that they do not differ much. This suggest following:
\begin{conjecture}
 The mixing time $\tau^{sep}_{mix}$ for \ctrt converges to  $n\log n$ as $n\to\infty$.
\end{conjecture}

\subsection{Optimal shuffling}
\ctrt is the shuffle used in \rcc and in Spritz. To reach stationarity (\ie 
produce random permutation), as we have shown, one needs to perform $O(n\log n)$ 
steps.
In each step we use a random number from interval $[0,\ldots,n-1]$, thus this shuffling requires  $O(n\log^2 n)$ random bits. \par 
\noindent 
One can ask the following question: Is this shuffling optimal in terms of required bits? The answer is no.
The entropy of the uniform distribution on $[n]$ is $O(n\log n)$ (since there are $n!$ permutations of $[n]$),
thus one could expect that optimal shuffling would require this number of bits.
Applying described \riffle  yields SST with  $ET=2\lg n$, at each step $n$ random bits are used, thus this shuffling requires $2n\lg n$ random bits,
matching the requirement of optimal shuffle (up to a constant).

In Figure \ref{fig:bit_compl} we summarized required number of bits for various 
algorithms and stopping rules.

\begin{figure}[H]
 \begin{center}
{\footnotesize
\begin{tabular}{l| c |  c | c | c}
 \# bits  & RC4 & Mironov & \ksa{${}_{KLZ}$} & RS \\
 \hline 
 used  & $40-2048$ &  &  &  \\
 \hline 
 \multirow{1}{*}{asympt.} &  & $(2 n H_n - 1)\lg n$ & $n 
(H_n + 1) \lg n$ & $2 n \lg n$ \\
 \hline 
 required   &  & $23~037$ & $14~590$ & $4~096$ \\
\end{tabular}
}
 \end{center}
\caption{Comparison between number of bits used by \rcc (40 to 2048) and 
required by mathematical models (Mironov~\cite{Mironov2002} and ours -- 
\ksa{${}_{KLZ}$}) versus 
length of the key for the time-reversed 
riffle shuffle. \textit{Bits asymptotics} approximates the number of fresh bits 
required by the mathematical model (number of bits required by the underlying 
Markov chain to converge to stationary distribution). \textit{Bits required} 
is (rounded) value of \textit{\# bits asymptotics} when $n = 256$.}\label{fig:bit_compl}
\end{figure}

\subsection{Perfect shuffle}\label{sec:perf_shuffle}
Let $\X$  be an ergodic chain with a stationary distribution $\psi$. Ergodicity means that nevertheless of initial distribution,
the distribution $\mathcal{L}(X_k)$ converges, in some sense, to $\psi$ as $k\to\infty$
(it is enough to think of $d_{TV}(\mathcal{L}(X_k),\psi)$ converging to 0 as $k\to\infty$).
Beside some trivial situations, the distribution of $\mathcal{L}(X_k)$ is \textsl{never} equal to 
$\psi$ for any finite $k$. That is why we need to study mixing time to know how close $\mathcal{L}(X_k)$ to $\psi$ is.
However, there are methods for so-called \textbf{perfect simulations}. This term refers to the art of converting a Markov chain
into the algorithm which (if stops) returns exact (unbiased) sample from stationary distribution of the chain.
The most famous one is the Coupling From The Past (CFTP) algorithm introduced in \cite{Propp1996}. The ingenious idea
is to realize the chain as evolving \textsl{from the past}. However, to apply algorithm effectively the chain must be
so-called \textsl{realizable monotone}, which is defined w.r.t some underlying partial ordering. Since we consider random walk
on permutations, applying CFTP seems to be a challenging task.
\par 
However, having Strong Stationary Time rule we can actually \textsl{perform perfect simulation}(!).
We simply run the chain until event SST occurs. At this time, the Definition \ref{def:SST} implies,
that the distribution of the chain is stationary. \par 
It is worth mentioning that all SSTs mentioned at the beginning of this Section are deterministic
(i.e., they are deterministic functions of the path of the chain, they do not use any extra randomness -
it is not the case for general SST). 
For a method for perfect simulation based on SST (and other methods) see \cite{LorekMarkowski2014}, where mainly monotonicity requirements
for three perfect sampling algorithms are considered.


\section{(Not so) random shuffles of RC4 -- revisited }\label{sec:ksa_proposed}
 
\medskip\par
\begin{figure}[h]
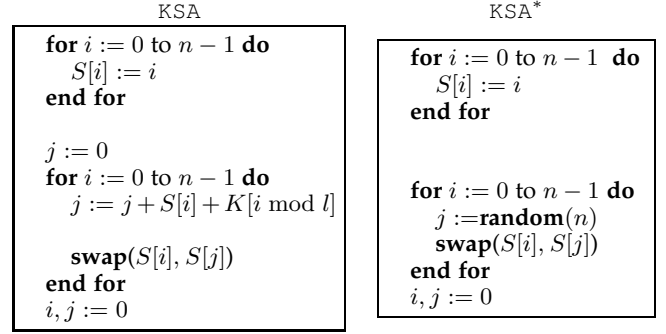

{\small
\begin{center}
\begin{tabular}{ccccc}
 \texttt{KSA}
&
\texttt{KSA${}^{*}$}
\\
\noindent\fbox{%
\begin{minipage}{.23\textwidth}
 \begin{algorithmic}
 \FOR{$i:=0$ to $n-1$}
\STATE $S[i]:=i$
\ENDFOR
\STATE \
\STATE $j:=0$
 \FOR{$i:=0$ to $n-1$}
 \STATE $j:=j+S[i]+K[i \bmod l]$
 \STATE \textbf{swap}($S[i],S[j]$)
 \ENDFOR
 \STATE $i,j:=0$
\end{algorithmic}
\end{minipage}%
}
 & 
\noindent\fbox{%
\begin{minipage}{.19\textwidth}
 \begin{algorithmic}
 \FOR{$i:=0$ to $n-1$ }
\STATE $S[i]:=i$
\ENDFOR
\STATE \
\STATE \ 
 \FOR{$i:=0$ to $n-1$}
 \STATE $j:=$\textbf{random}$(n)$
 \STATE \textbf{swap}($S[i],S[j]$)
 \ENDFOR
 \STATE $i,j:=0$
\end{algorithmic}
\end{minipage}%
}
\end{tabular}
\end{center}
\caption{\ksa of \rcc~ algorithm and its idealized version $\ksa^*$.
}
\label{fig:KSA_and_idealizedb}
}
\end{figure}

\noindent 
\textbf{\rcc algorithm.}
\rcc is a stream cipher, its so-called \textsl{internal state} is $(S,i,j)$, 
where $S$ is a permutation of $[n]$ and 
$i,j$ are some two indices. As input it takes $L-$byte message $m_1, \ldots, 
m_L$ and a secret key $K$ and returns 
ciphertext $c_1,\ldots,c_L$. The initial state is the output of \ksa. Based on this state
\prga   is used to output bits which are 
XORed with the message. 
The actual \ksa algorithm used in \rcc is presented in  Figure \ref{fig:KSA_and_idealizedb} together with its \textsl{idealized version} \ksa{${}^{*}$}
(where a secret key-based randomness is replaced with pure randomness)

The goal of \ksa of original \rcc is to produce a pseudorandom permutation of $n=256$ cards.
The original   algorithm performs $256$ steps of \ctrt. 
Even if we replace $j:=j+S[i]+K[i \bmod l]$ in \ksa with $j=\textbf{random}(n)$ (i.e., 
we actually consider idealized version), then  $\Theta(n\log n)$ steps are needed.
Considering idealized version was a main idea of Mironov \cite{Mironov2002}. 
Moreover, he showed that $2n\log n$ steps are enough. 
However, our result stated in Lemma \ref{lem:ksa_ET} together 
with experimental results from Section \ref{sec:expResults} imply that around $n\log n$ steps are enough.
Moreover, if we consider \ksa with  \rtrt (instead of \ctrt) then we can state more formal theorem in terms of  
introduced security definitions.
\begin{theorem}\label{twr:rc4_main}
Let $\mathcal{A}$ be an adversary. Let $K\in\{0,1\}^{r n}$ be a secret key. Let $\mathcal{S}(K)$ be 
\ksa with \rtrt shuffling  which runs for 
$$ r= n(H_n+1)+{\pi n\over 2} {1\over \sqrt{n!\varepsilon}}$$
steps with $0<\varepsilon<{1\over n!}$. Then 
\[\left|
\underset{\scriptscriptstyle K \leftarrow \{0, 1\}^{rm}}{Pr}
 \left[\mathcal{A}( \textsf{S}(K)) = 1 \right] -  
{\underset{\scriptscriptstyle R \leftarrow 
\mathcal{U}(\mathcal{S}_n)}{Pr}}[\mathcal{A}( 
R) = 1] \right| \leq \varepsilon. 
\]

\end{theorem}

\noindent 
\textsl{Proof}.
The stopping rule \texttt{StoppingRuleKLZ} given in Section \ref{sec:betterSST} is an SST for \rtrt, denote it by $T$. 
From  Lemma \ref{lem:ksa_ET} we know that $ET=n(H_n+1)+O(H_n)$ and $VarT\sim {\pi^2\over 4}n^2$. 
Inequality (\ref{eq:TV_sep_SST}) and Chebyshev's inequality imply that for $r=n(H_n+1)+{\pi n\over 2} c$ we have 
$$sep(\mathcal{L}(X_r),\psi)\leq P(T>r)\leq {1\over c^2},$$
\ie $r=\tau_{mix}^{sep}(c^{-2})$.
Taking $c={1\over \sqrt{n!\varepsilon}}$ we perform   \rtrt for $r=\tau_{mix}^{sep}(n!\varepsilon)$ steps, \ie $sep(\mathcal{L}(X_r),\psi)\leq n!\varepsilon$ 
Inequality (\ref{eq:sep_unif}) implies that $|Pr(X_r=\sigma)-{1\over n!}|\leq \varepsilon$ for any permutation $\sigma$ and thus completes the proof.
\QEDB

\subsection{Sign distinguisher for or \ctrt.}\label{sec:signDist}
\textbf{Sign distinguisher.} For a permutation $\pi$ which has a 
representation of non-trivial transpositions $\pi = (a_1 b_1)(a_2 b_2) \ldots 
(a_m b_m)$ the sign is defined as:
$sign(\pi) = (-1)^m$. So the value of the sign is $+1$ whenever $m$ is even and 
is equal to $-1$ whenever $m$ is odd. The distinguisher is presented in Figure \ref{fig:rc4_sign_dist}.

\begin{center}
\begin{figure}[H]\caption{Sign distinguisher for \ctrt}
\begin{center}
\fbox{ 
\begin{varwidth}{\dimexpr\linewidth-2\fboxsep-2\fboxrule}
\begin{algorithmic}
 \STATE Input: $S, t$
 \STATE Output: $b$
 \STATE \textbf{if} sign($S$)=$(-1)^t$ 
 \STATE \ \ \textbf{then return} \textit{false} 
 \STATE \ \ \textbf{else return} \textit{true}
\end{algorithmic}
\end{varwidth}%
}
\end{center}\label{fig:rc4_sign_dist}
\end{figure}
\end{center}

It was observed in~\cite{Mironov2002} that the sign of the 
permutation at the 
end of \rcc's \ksa algorithm is not uniform. 
And as a conclusion it was noticed that the number of discarded shuffles (by 
\prga) must grow at 
least linearly in $n$. 
Below we present this result obtained in a different way than 
 in~\cite{Mironov2002}, giving the exact formula for advantage 
at any step $t$.


One can look at the sign-change process for the \ctrt~ as follows: after the 
table is 
initialized,  sign of the permutation is $+1$ since it is 
identity so the initial distribution is concentrated in $v_0 = 
(Pr(\textrm{sign}(Z_0)=+1),Pr(\textrm{sign}(Z_0)=-1))=(1, 0)$. 

Then in each step the sign is unchanged if and only if  $i = j$ 
which happens 
with probability $1/n$. So the transition matrix $M_n$ of a sign-change 
process induced by the shuffling  process is equal to:
\[
      M_n := \left(
        \begin{array}{cc}
         \frac{1}{n} & 1-\frac{1}{n}\\ 
         1-\frac{1}{n} & \frac{1}{n}
         \end{array}
      \right).
\]
This conclusion corresponds to 
looking at the distribution of the sign-change process after $t$ steps: $v_0 
\cdot M_n^t 
$, where $v_0$ is the initial distribution.
The eigenvalues and eigenvectors of $M_n$ are $(1,{2-n\over n})$ and $(1,1)^T,(-1,1)^T$ respectively. The spectral decomposition yields

$$\scriptsize
\begin{array}{llll}
v_0\cdot M_n^t  & = & (1,0)
\left(
\begin{array}{rr}
 1 & -1 \\[8pt]
 1 & 1
\end{array}
\right)
\left(
\begin{array}{rr}
 1 & 0 \\[8pt]
 0 & {2-n\over n}
\end{array}
\right)^t
\left(
\begin{array}{rr}
 {1\over 2} & -{1\over 2} \\[8pt]
 -{1\over 2} & {1\over 2}
\end{array}
\right) \\[15pt]
 & = &
 \displaystyle \left({1\over 2}+{1\over 2}\left({2\over n}-1\right)^t,{1\over 2}-{1\over 2}\left({2\over n}-1\right)^t \right).
\end{array}
$$
In Fig.~\ref{tab:signAdv} 
the advantage $\epsilon$ of a 
sign-adversary after dropping $k$ bytes of the output (so after $n + k$ steps 
of the shuffle, for the mathematical model) is presented.
\begin{figure}[H]
\begin{center}
\begin{tabular}[c]{r | c | c | r}
$k$ & $+1$ & $-1$ & $\epsilon$\\
\hline 
0 & .5671382998250798 & .4328617001749202 & $2^{-3.89672}$\\
256 & .509015 & .490985 & $2^{-6.79344}$\\
512 & .5012105173235390 & .4987894826764610 & $2^{-9.69016}$\\
768 & .500163 & .499837 & $2^{-12.5869}$\\
1024 & .5000218258757580 & .4999781741242420 & $2^{-15.4836}$\\
2048 & .5000000070953368 & .4999999929046632 & $2^{-27.0705}$\\
4096 & .5000000000000007 & .4999999999999993 & $2^{-50.2442}$\\
8192 & .5 & .5 & $2^{-96.5918}$\\
\end{tabular}
\end{center}
\caption{The advantage ($\epsilon$) of the Sign distinguisher of \rcc 
after $n 
+ k$ steps or equivalently, after  
discarding initial $k$ bytes.}
\label{tab:signAdv}
\end{figure}
For $n = 256$ (which corresponds to the value of $n$ used in \rcc)
and initial 
distribution being identity permutation after $t = n = 256$ steps one gets:
$v_0 \cdot M_{256}^{256} = (0.567138, 0.432862)$.

In~\cite{ManSha2001} it was suggested that the first 512 bytes of the output should 
be dropped.

%


\section{Strong stationary time based \ksa algorithms}\label{sec:nasz_ksa}
In addition to \ksa and $\ksa^*$ (given in Fig. \ref{fig:KSA_and_idealizedb}) we  introduce another version called \ksa{${}^{**}_\texttt{Shuffle,ST}$($n$)}.
The main difference is that it does not run a predefined number of steps, but it takes some procedure $\texttt{ST}$ as a parameter
(as well as \texttt{Shuffle} procedure). The idea is that it must be a rule corresponding to Strong Stationary Time for given shuffling scheme.
The  \ksa{${}^{**}_\texttt{Shuffle,ST}$($n$)} is given in the following algorithm:
\begin{algorithm}[H]
\caption{\ksa{${}^{**}_\texttt{Shuffle,ST}$($n$)}}\label{alg:ksa2}
\begin{algorithmic}
\REQUIRE Card shuffling \texttt{Shuffle} procedure, stopping rule \texttt{ST} which is a Strong Stationary Time for \texttt{Shuffle}.
 \STATE 
 \FOR{$i:=0$ to $n-1$}
\STATE $S[i]:=i$
\ENDFOR
\STATE \
 \WHILE{($\neg $ \texttt{ST})}
 \STATE \texttt{Shuffle}($S$)
 \ENDWHILE
\end{algorithmic}
\end{algorithm}
The  algorithm   works as follows. It starts with identity permutation. Then at 
each step it performs 
some card shuffling procedure \texttt{Shuffle}. 
Instead of running it for a 
pre-defined number of steps, it runs until an event defined by a procedure 
\texttt{ST} occurs. The procedure \texttt{ST} is designed in such a 
way that it guarantees that the event is a Strong 
Stationary Time. At each step the algorithm uses new randomness -- one can 
think 
about that as of an idealized version but when the length of a key is greater than 
the number of random bits required by the algorithm then we end up with a 
permutation which cannot be distinguished from a random (uniform) one (even by 
a computationally unbounded adversary).

\noindent
\textbf{Notational convention}: in  \ksa{${}^{**}_\texttt{Shuffle,ST}$($n$)} we omit parameter $n$. Moreover,
if \texttt{Shuffle} and \texttt{ST} are omitted it means that we use \ctrt as shuffling procedure and stopping rule is clear 
from the context.
Note that if we use for stopping rule  \texttt{ST} ``stop after $n$ steps'' 
(which of course is \textsl{not} SST), it is equivalent to \rcc's 
\ksa{${}^{*}$}.

Given a shuffling procedure \texttt{Shuffle} one wants to have a ``fast'' 
stopping rule \texttt{ST}
(perfectly one wants an optimal SST which is stochastically the smallest).
The stopping rule \texttt{ST} is a parameter, since for a given shuffling 
scheme one can come up with a better stopping rule(s).
This is exactly the case with \ctrt and \rtrt. In Section \ref{sec:shuffling} we recalled Mironov's 
\cite{Mironov2002} stopping rule as well as new introduced ``faster''
rule called (called \texttt{StoppingRuleKLZ}).
\medskip \par \noindent
Let us point out few important properties of \ksa{${}^{**}_\texttt{Shuffle,ST}$($n$)} algorithm:\medskip\par
\begin{itemize}
 \item[\textbf{P1}] The resulting permutation is purely random. \medskip\par
 \item[\textbf{P2}] The information about the number of rounds that 
an SST-algorithm performs does not reveal any information about the resulting 
permutation.
\end{itemize}
The property \textbf{P1} is due to the fact that we actually perform \textsl{perfect simulation} (see Section \ref{sec:perf_shuffle}).
The property \textbf{P2} follows from the fact that the Definition \ref{def:SSTa} is equivalent to the following one 
(see \cite{Aldous1987}):
\begin{definition}\label{def:SSTa}
Random variable $T$ is \textbf{Strong Stationary Time (SST)} if it is a
randomized stopping time for chain $\X$ such that:
$$
X_T \textrm{ has distribution } \psi \textrm{ and is independent of } T.
$$ 
 \end{definition}
 
The property \textbf{P1} is the desired property of the output of the \ksa algorithm,
whereas \textbf{P2} turns out to be crucial for our timing attacks applications (details in Section \ref{sec:timing_sst}).


\section{Timing attacks and SST}\label{sec:timing_sst}
 In the black box model (for CPA-security experiment), an adversary has access 
to an oracle computing encryptions $F(k,\cdot)$ for a secret key $k$. So the 
adversary is allowed to obtain a vector of cryptograms ($c_1, \ldots, c_q$) for 
plaintexts ($m_1, \ldots, m_q$) of its choice (where each $c_i = \tuple{r_i, 
y_i} = F(k, m_i, r_i)$). Each $r_i$ is the randomness used by the oracle to 
compute encryption (\eg its initial counter \textit{ctr} for CTR-mode, 
\textit{IV} for CBC-mode \etc).
So the adversary has a vector of tuples 
\[\mathcal{O}_{\text{BlackBox}, F} = \tuple{m_i, c_i} = \tuple{m_i, F(k, m_i, 
r_i)}.
\]

In the real world, as a result of interacting with an implementation of 
$F$, a timing information may also be accessible to an adversary. 
In such a case her knowledge is a vector of tuples 
\[ 
\mathcal{O}_{Timing, F} = \tuple{m_i, c_i, t_i}  = 
\tuple{m_i, F(k, 
m_i, r_i), t_i},
\]
where $t_i$ is the time to reply with $c_i$ on query $m_i$ and 
randomness used $r_i$.
(A very similar information is accessible to an adversary in 
chosen-ciphertext attacks.)

An encryption scheme is susceptible to a timing attack when there exists an 
efficient function $H$ which on input $H\left(\tuple{m_1, c_1, t_1}, \ldots, 
\tuple{m_w, c_w, t_w}\right)$ outputs $\lambda(k)$ (some function of a secret 
key $k$). Depending on the scheme and on a particular implementation the 
leakage function can be of various forms \eg $\lambda(k) = HW(k)$ -- the 
Hamming weight of $k$, or in the best 
 case $\lambda(k) = k$.

The reason why a timing attack on implementation of a function $F$  ($y = 
F(k, x)$) is possible is because the running time of $F$ is not the same on 
each input.
Usually timing attacks focus on a particular part of $F$ -- $f_{weak}$ which 
operates on a function $h_x$ of the original input $x$ and on a function $h_k$ 
of the original key (\eg $h_k$ may be a substring of a round key). When 
$F(k, x) = f_{rest}(k, f_{weak}(k, x))$ the statistical model is built and 
$T(F(k, x)) = T(f_{weak}(k, x)) + T(f_{rest}(k, f_{weak}(k, x)))$. The 
 distribution of $T(f_{weak}(k, x))$ is conditioned on $k$ and estimating this 
value (conditioned on subset of bits of $k$) is the goal of the adversary.

\subsection{Countermeasures -- approach I}
There is one common approach to eliminate possibility of timing attacks, namely 
to re-implement $F$ ($ct(F) = F_{ct}$ -- constant-time implementation of $F$) 
in such a way that for all pairs $(k, x) \neq (k', x')$ the 
running times are the same $T(F(k, x)) = T(F(k', x')) = t$.
Then, an adversary's observation:
\[ 
\mathcal{O}_{Timing, ct(F)} = \tuple{m_i, c_i, t}  = 
\tuple{m_i, F_{ct}(k, 
m_i, r_i), t}
\]
becomes the same as in the black-box model.

\subsection{Countermeasures -- approach II}

The other approach is to still use a weak implementation, but ``blind'' (mask) 
the 
adversary's view. Let us start with a public-key example, described in the 
Section~\ref{sec:rw-masking}. The RSA's fragile decryption $F(d, c) = 
f_{weak}(d, c)$ (Algorithm~\ref{alg:modpower}) is implemented with masking 
(Algorithm~\ref{alg:rsasecureG}) as:
$F_{masked}(d, c) = g_2(f_{weak}(d, g_1(c, r)), r)$, where $g_1(c, r) = c r^e 
\bmod N$ and $g_2(x, r) = x r^{-1} \bmod N$.
In this case, the adversary's view is:
\[ 
\mathcal{O}_{Timing, masked(F)} =  
\tuple{c_i, F_{masked}(d, 
c_i), t_i}.
\]
While $t_i$ are again different, they depend on both $d$ and $g_1(c_i, r_i)$,
but an adversary  does not obtain information about random $r_i$'s and thous 
information about $t_i$ is useless.

\subsection{Non-constant time implementation immune to timing attacks for 
symmetric encryption}

Let us use the similar approach for the symmetric 
encryption case.  For a time-attack susceptible function $F$ we would like to 
build a function $G$ by designing functions $g_1(\cdot)$ and $g_2(\cdot)$ and 
then to define 
\[
G(k, x) = g_2(F(k, g_1(x, \cdot)), \cdot).
\] 
We restrict 
ourselves only to functions $G$ of that form. We do not require that $G(k, x) = 
F(k, x)$ as it was in the public key case. It would be enough for us that there 
exists an efficient algorithm for decryption -- this is satisfied if both 
$g_1, g_2$ have efficiently computable inversions.

Now let us consider the following  
approaches for selecting $g_1$ (and $g_2$):
\begin{enumerate}
 \item $g_1(x)$ -- then such an approach is trivially broken because an 
attacker knows exactly what is the input to $F$.

 \item $g_1(x, r)$ for a randomly selected $r$ generated during run of $G$. 
Assuming that $r$ cannot be efficiently removed from $F(k, g_1(x, r))$ via 
application of $g_2$, the ciphertext of $G$ would need to include $r$ but then 
this is exactly the same case as (1). 
 
 \item $g_1(k, x)$ -- in that case an attacker might want to perform a 
timing attack on $g_1$ -- one would need to guarantee that $g_1$~itself is 
immune to timing attacks. This approach is used in McBits~\cite{McBits} 
system, $g_1$~is a constant-time implementation of Bene\v{s} permutation 
network \cite{BenesNetwork}.
 
 \item $g_1(k, x, r)$ -- again value $r$ would need to be included in 
ciphertext generated by $G$ but now for an attacker this case would be no 
easier than when $g_1$ depends solely on $k$ (just as in the case 3).
 
\end{enumerate}

While $g_2$ cannot remove randomness (unless $F$ is 
trivially broken) and the randomness used in $g_1$ needs to be included in the 
ciphertext, one needs $g_2$ anyway. It is required for the protection in the 
case when attacker has access to the decryption oracle -- then $g_2$ plays the 
same role as $g_1$ plays for protecting against timing attacks during 
encryption.

We suggest the following approach, let $k = (k_f, k_g)$ where $k_f$ is a 
portion of the secret key used by $F$:
\begin{enumerate}
 \item Let $g_1(k_g, x), g_2(k_g, x)$ be permutations of $n=|x|$ bits depending on a key. \textsl{I.e.,} 
 $g_1(k_g,\cdot)$ and $ g_2(k_f,\cdot)$ are permutations constructed based on $k_g$ and $k_f$ respectively,
 which are later on applied to $x$, a sequence of $n$ bits (similarly as in McBits scheme)
 \item Both $g_1, g_2$ are the permutations obtained as a result of 
\ksa{${}^{**}_\texttt{Shuffle,ST}$} algorithm -- one has a freedom in the 
choice of a shuffle algorithm.

\end{enumerate}

\begin{figure}[H]
\begin{center}
\tikzstyle{decision} = [diamond, draw, fill=white!100, 
    text width=4.5em, text badly centered, node distance=1cm, inner sep=0pt]
\tikzstyle{block} = [rectangle, draw, fill=white!100,
    text width=7em, text centered, rounded corners, minimum height=2em]
\tikzstyle{block2} = [rectangle, draw, fill=white!100, 
    text width=5em, text centered, rounded corners, minimum height=4em]    
\tikzstyle{line} = [draw, -latex']
\tikzstyle{cloud} = [draw, ellipse,fill=red!20, node distance=1cm,
    minimum height=2em]
    
\begin{tikzpicture}[node distance = 1.4cm, auto]
    \node  (init) {$x$};    
    \node [block, below of=init] (comp) {$g_1$};
    \node [left of=comp, node distance=2.4cm] (key) {$k_g$};
    \node [block, below of=comp] (F_bad) {\small $F$};
    \node [left of=F_bad, node distance=2.4cm] (key2) {$k_f$};
    \node [block, below of=F_bad] (comp2) {$g_2$};
    \node [left of=comp2, node distance=2.4cm] (key3) {$k_g$};
    
    \node [below of=comp2] (output) {\small $y$};
%
%
      \path [line] (key) -- (comp);
      \path [line] (init) -- (comp);
      \path [line] (comp) -- (F_bad);       
      \path [line] (key2) -- (F_bad);       
      \path [line]  (F_bad) -- (comp2);       
      \path [line]  (key3) -- (comp2);       
      \path [line]  (comp2) -- (output);  

\end{tikzpicture}
\caption{Timing attack immune $G$ obtained from timing attack 
susceptible $F$ via compounding with SST-based shuffling algorithms $g_1, 
g_2$.} \label{fig:timing-immune}
\end{center}
\end{figure}
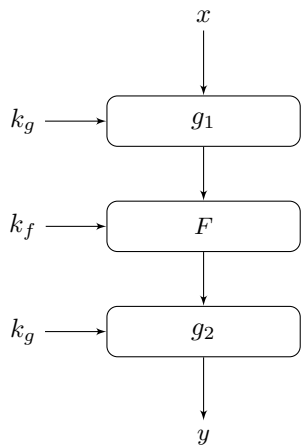

\begin{algorithm}[H]
\caption{$G$ -- timing attack immune encryption based on 
timing attack susceptible symmetric encryption $F$}
\label{alg:rsasecureG}
\begin{algorithmic}[1]
\REQUIRE $x, F, g_1, g_2, k_f, k_g$
\STATE compute $y_1 = g_1(k_g, x)$
\STATE compute $y_2 = F(k_f, y_1)$
\RETURN $y = g_2(k_g, y_2)$
\end{algorithmic}
\end{algorithm}

\begin{algorithm}[H]
\caption{$G^{-1}$ -- timing attack immune decryption implementation based on 
timing attack susceptible symmetric decryption $F^{-1}$}
\label{alg:rsasecureG1}
\begin{algorithmic}[1]
\REQUIRE $y, F^{-1}, g^{-1}_1, g^{-1}_2, k_f, k_g$
\STATE compute $x_1 = g^{-1}_2(k_g, y)$
\STATE compute $x_2 = F^{-1}(k_f, x_1)$
\RETURN $x = g^{-1}_1(k_g, x_2)$
\end{algorithmic}
\end{algorithm}
\noindent
Resulting $G$ has the following properties:
\begin{description}
 \item[\textbf{G1}] Resulting permutation $g_1$ is uniformly random, provided a long 
enough 
$k_g$ is applied (see property \textbf{P1} from Section~\ref{sec:nasz_ksa}).
 \item[\textbf{G2}] For a given $k_g$ and for any $x, x'$ one obtains $T(g_1(k_g, x)) = 
T(g_1(k_g, x'))$ (and the same holds for $g_2$) -- while $g_1$ is a result of 
SST, running it for the same input $k_g$ leads to the same running time.
 \item[\textbf{G3}] The running time of $g_1$ (and of $g_2$) does not reveal any 
information about the 
resulting permutation (see property \textbf{P2}). As a result, the 
adversary does not have influence on inputs to~$F$.
\end{description}

The property \textbf{G3} guarantees security of the outputs of $g_1$ and $g_2$ 
and that they cannot be derived from the running time but it does no guarantee 
that $k_g$ will remain secure. It turns out  that some bits of $k_g$ may leak 
through the information about the running time but:
\begin{itemize}
 \item there are at least $n!$ (where $n$ is the number of bits of $x$) different keys which have the same running time,
 \item the information which leaks from the running time does not let to obtain 
ANY information about the resulting permutation.
\end{itemize}
For more details see 
Section~\ref{app:timing}.

\subsection{Implementation}
To implement the above construction, one needs (pseudo)random bits to generate 
a permutation $g_1$. These bits can be obtained in the similar way one may 
produce a key for authenticated encryption: 
let $k$ be a 128-bit AES key, calculate $k_g \leftarrow$ {AES}$(k,0)$, 
 $k_g \leftarrow$ {AES}$(k,1)$, and $k_m \leftarrow$ {AES}$(k,2)$. Then 
{AES}$(k_g, IV)$ is used as a stream cipher to feed \riffle until SST occurs, 
producing $g_1$. The following bits in the same 
way produce $g_2$.  $F$ is \eg {AES}-128 in CBC or CTR mode, all having 
$k_f$ as a secret key. Then, $k_m$ may be used to compute the MAC (this 
disallows an adversary to play with \eg hamming weight of inputs of $g_1$).


A proof of concept implementation was prepared (as a student project: 
\url{https://bitbucket.org/dmirecki/bearssl/}). This implementation first 
uses the AES key and riffle-shuffle to generate an array 
(\textsf{permutation}) which stores the permutation $g_1$ (of size 128).
Then it permutes bits of the buffer by performing 128-step loop (line 34 of the 
listing~\ref{code:riffleSlow}).

\begin{listing}
\begin{minted}[label=/src/riffle.c,
               mathescape,
               linenos,
               firstnumber=34,
               fontsize=\footnotesize,
               numbersep=5pt,
               frame=lines,
               xleftmargin=15pt,
               framesep=2mm]{c}
for (uint32_t i = 0; i < 128; i++) {
   if (CHECK_BIT(buf[i >> 3], i & 7)) {
            SET_BIT(temp_buf[permutation[i] >> 3], 
	      permutation[i] & 7);
   }
}
\end{minted}
\caption{A fragment of the 
(\textsf{src/riffle.c}) responsible for applying a 
shuffle.}\label{code:riffleSlow}
\end{listing}
Fix $n=128=2^7$. Any permutation of $[n]$ can be represented as 7 steps of \textsl{Riffle Shuffle}.
This can be exploited to speed up the application of the permutation:
\begin{enumerate}
 \item generate a permutation using \textsl{Riffle Shuffle}  with SST (14 steps on average), say permutation $\sigma$
 was obtained.
 \item represent $\sigma$ as 7 steps of some \textsl{other} \textsl{Riffle Shuffle} scheme, \ie compute the corresponding 
 $7\cdot 128$ bits.
 Do the same for $\sigma^{-1}$.
 \item apply this  7-step \textsl{Riffle Shuffle} on  input bits.
\end{enumerate}

Then the last step may be implemented using efficient instructions: 
\texttt{\_pext\_u64, \_pdep\_u64, \_mm\_popcnt\_u64} or its 32-bit analogs.
 Consult \cite{Intel_perm} p. 4-270 -- 4-273 Vol. 2B for above instructions.

\begin{figure}[H]
\begin{center}
{\small
 \begin{tabular}[c]{l | c |  c | c | c}
  & big & ct & student & sst \\
  \hline
  AES-128 CBC encrypt & 170.78 & 29.06 & 11.46 & 37.63\\
  \hline
  AES-128 CBC decrypt & 185.23 & 44.28 & 11.48 & 37.65\\
  \hline
  AES-128 CTR & 180.94 & 55.60 & 16.83 & 37.52 \\
 \end{tabular}
}
\caption{Comparison of the efficiency of selected algorithms: \texttt{big} -- 
classic, table-based AES implementation, \texttt{ct} -- constant time 
implementation using bitslicing strategy, \texttt{student} --  
implementation of 
the sst-masking described in Section~\ref{sec:timing_sst} with slow 
permutation application, \texttt{sst} -- estimated implementation which of 
the sst-masking which uses BMI2 and SSE4.2 instructions. All results are in 
MB/s and were obtained by running speed tests of the 
BearSSL~\cite{Pornin} (via the command: \texttt{testspeed all} on Intel 
i7-4712HQ with 16GB RAM).}\label{tab:efficiency}
\end{center}
\end{figure}


\subsection{Timing attacks and \ksa{**}}\label{app:timing}

This section explains  why bits of the secret key 
$k_g$ used to generate $g_1, g_2$ need to be distinct to the bits $k_f$ which 
are used by $F$. This is a simple conclusion of properties 
\textbf{P1} and \textbf{P2} (defined in Section~\ref{sec:nasz_ksa}).
The examples are  intended to give a reader a better insights on the 
properties achieved by SST shuffling algorithms.
As already mentioned - regardless the duration of SST algorithm adversary can learn something 
about the key, but has completely no knowledge about the resulting permutation. 
The examples will also show how a choice of a shuffling algorithm influences the leaking information about the key.
In both cases we consider extreme (and very unlikely) situation where algorithm 
runs for the smallest number of steps.

\subsubsection{Top-to-random card shuffling}

Consider the algorithm  \ksa{${}^{**}_\texttt{T2R,ST}$} with the shuffling 
procedure 
corresponding to \textsl{Top-To-Random} card shuffling:
put the card which is currently on the top to the position $j$ defined by 
the  bits of the key $k_g$ in the current round.  Recall  stopping rule
\texttt{ST} which is SST for this shuffle: stop one step after the card which was initially at the bottom reaches the top of the 
deck.
   Before the start of the algorithm, mark the last card (\ie the card 
$n$ is marked\footnote{It is known \cite{Aldous1986} that  optimal \sst for 
\textit{top to random} initially marks card which is second from  the bottom.}).

\tikzstyle{RectObject}=[rectangle,fill=white,draw,line width=0.5mm]
\tikzstyle{line}=[draw]
\tikzstyle{arrow}=[draw, -latex]

\tikzstyle{block} = [draw,fill=gray!40,minimum size=0.5em]
\tikzstyle{ble} = [draw,minimum size=0.5em]
\tikzstyle{ble2} = [draw,fill=red!20,minimum size=0.5em]

\tikzstyle{block2} = [draw,fill=gray!90,minimum size=0.5em]
\def\radius{.7mm} 
\tikzstyle{branch}=[fill,shape=circle,minimum size=3pt,inner sep=0pt]
 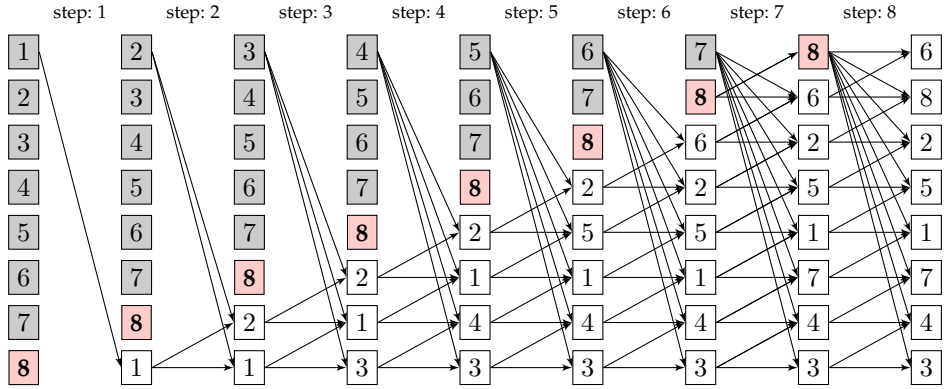
\begin{figure*}
\begin{center}
\begin{tikzpicture}[>=latex',scale=1.0,every node/.style={scale=1}]

   \node  at (-.25,-.5) {\scriptsize step: 1};
   \node  at (1.25,-0.5) {\scriptsize step: 2};
   \node  at (2.75,-0.5) {\scriptsize step: 3};
   \node  at (4.25,-0.5) {\scriptsize step: 4};
   \node  at (5.75,-0.5) {\scriptsize step: 5};
   \node  at (7.25,-0.5) {\scriptsize step: 6};
   \node  at (8.75,-0.5) {\scriptsize step: 7};
   \node  at (10.25,-0.5) {\scriptsize step: 8};

  \node[block] at (-1,-1)   (ak1) {$1$};
  \node[block] at (-1,-1.6) (ak2) {$2$};
  \node[block] at (-1,-2.2) (ak3) {$3$};
  \node[block] at (-1,-2.8) (ak4) {$4$};
  \node[block] at (-1,-3.4) (ak5) {$5$};
  \node[block] at (-1,-4)   (ak6) {$6$};
  \node[block] at (-1,-4.6)   (ak7) {$7$};
  \node[ble2] at (-1,-5.2)   (ak8) {$\textbf{8}$};

 \node[block] at (.5,-1)    (bk1) {$2$};
  \node[block] at (.5,-1.6) (bk2) {$3$};
  \node[block] at (.5,-2.2) (bk3) {$4$};
  \node[block] at (.5,-2.8) (bk4) {$5$};
  \node[block] at (.5,-3.4) (bk5) {$6$};
  \node[block] at (.5,-4)   (bk6) {$7$};
  \node[ble2] at (.5,-4.6)   (bk7) {$\textbf{8}$};
  \node[ble] at (.5,-5.2)   (bk8) {$1$};
   \draw[->] (ak1.east) to   (bk8.west) ;

 \node[block] at (2,-1)    (ck1) {$3$};
  \node[block] at (2,-1.6) (ck2) {$4$};
  \node[block] at (2,-2.2) (ck3) {$5$};
  \node[block] at (2,-2.8) (ck4) {$6$};
  \node[block] at (2,-3.4) (ck5) {$7$};
  \node[ble2] at (2,-4)   (ck6) {$\textbf{8}$};
  \node[ble] at (2,-4.6)   (ck7) {$2$};
  \node[ble] at (2,-5.2)   (ck8) {$1$};

  \draw[->] (bk1.east) to   (ck7.west) ;
  \draw[->] (bk1.east) to   (ck8.west) ;
   
  \draw[->] (bk8.east) to   (ck7.west) ;
  \draw[->] (bk8.east) to   (ck8.west) ;
 \node[block] at (3.5,-1)    (dk1) {$4$};
  \node[block] at (3.5,-1.6) (dk2) {$5$};
  \node[block] at (3.5,-2.2) (dk3) {$6$};
  \node[block] at (3.5,-2.8) (dk4) {$7$};
  \node[ble2] at (3.5,-3.4) (dk5) {$\textbf{8}$};
  \node[ble] at (3.5,-4)   (dk6) {$2$};
  \node[ble] at (3.5,-4.6)   (dk7) {$1$};
  \node[ble] at (3.5,-5.2)   (dk8) {$3$};

  \draw[->] (ck1.east) to   (dk6.west) ;
  \draw[->] (ck1.east) to   (dk7.west) ;
  \draw[->] (ck1.east) to   (dk8.west) ;

  \draw[->] (ck7.east) to   (dk6.west) ;
  \draw[->] (ck7.east) to   (dk7.west) ;
   
  \draw[->] (ck8.east) to   (dk7.west) ;
  \draw[->] (ck8.east) to   (dk8.west) ;

 \node[block] at (5,-1)    (ek1) {$5$};
  \node[block] at (5,-1.6) (ek2) {$6$};
  \node[block] at (5,-2.2) (ek3) {$7$};
  \node[ble2] at (5,-2.8) (ek4) {$\textbf{8}$};
  \node[ble] at (5,-3.4) (ek5) {$2$};
  \node[ble] at (5,-4)   (ek6) {$1$};
  \node[ble] at (5,-4.6)   (ek7) {$4$};
  \node[ble] at (5,-5.2)   (ek8) {$3$};

  \draw[->] (dk1.east) to   (ek8.west) ;
  \draw[->] (dk1.east) to   (ek7.west) ;
  \draw[->] (dk1.east) to   (ek6.west) ;
  \draw[->] (dk1.east) to   (ek5.west) ;

  \draw[->] (dk6.east) to   (ek5.west) ;
  \draw[->] (dk7.east) to   (ek6.west) ;
  \draw[->] (dk8.east) to   (ek7.west) ;

  \draw[->] (dk6.east) to   (ek6.west) ;
  \draw[->] (dk7.east) to   (ek7.west) ;
  \draw[->] (dk8.east) to   (ek8.west) ;

 \node[block] at (6.5,-1)    (fk1) {$6$};
  \node[block] at (6.5,-1.6) (fk2) {$7$};
  \node[ble2] at (6.5,-2.2) (fk3) {$\textbf{8}$};
  \node[ble] at (6.5,-2.8) (fk4) {$2$};
  \node[ble] at (6.5,-3.4) (fk5) {$5$};
  \node[ble] at (6.5,-4)   (fk6) {$1$};
  \node[ble] at (6.5,-4.6)   (fk7) {$4$};
  \node[ble] at (6.5,-5.2)   (fk8) {$3$};

  \draw[->] (ek1.east) to   (fk4.west) ;
  \draw[->] (ek1.east) to   (fk5.west) ;
  \draw[->] (ek1.east) to   (fk6.west) ;
  \draw[->] (ek1.east) to   (fk7.west) ;
  \draw[->] (ek1.east) to   (fk8.west) ;

  \draw[->] (ek8.east) to   (fk7.west) ;
  \draw[->] (ek8.east) to   (fk8.west) ;

    \draw[->] (ek7.east) to   (fk6.west) ;
    \draw[->] (ek6.east) to   (fk5.west) ;
    \draw[->] (ek5.east) to   (fk4.west) ;
    \draw[->] (ek7.east) to   (fk7.west) ;
    \draw[->] (ek6.east) to   (fk6.west) ;
    \draw[->] (ek5.east) to   (fk5.west) ;

 \node[block] at (8,-1)    (gk1) {$7$};
  \node[ble2] at (8,-1.6) (gk2) {$\textbf{8}$};
  \node[ble] at (8,-2.2) (gk3) {$6$};
  \node[ble] at (8,-2.8) (gk4) {$2$};
  \node[ble] at (8,-3.4) (gk5) {$5$};
  \node[ble] at (8,-4)   (gk6) {$1$};
  \node[ble] at (8,-4.6)   (gk7) {$4$};
  \node[ble] at (8,-5.2)   (gk8) {$3$};


  \draw[->] (fk1.east) to   (gk3.west) ;
  \draw[->] (fk1.east) to   (gk4.west) ;
  \draw[->] (fk1.east) to   (gk5.west) ;
  \draw[->] (fk1.east) to   (gk6.west) ;
  \draw[->] (fk1.east) to   (gk7.west) ;
  \draw[->] (fk1.east) to   (gk8.west) ;
   
  \draw[->] (fk7.east) to   (gk6.west) ;
  \draw[->] (fk7.east) to   (gk7.west) ;

  \draw[->] (fk8.east) to   (gk7.west) ;
  \draw[->] (fk8.east) to   (gk8.west) ;

  \draw[->] (fk4.east) to   (gk3.west) ;
  \draw[->] (fk5.east) to   (gk4.west) ;
  \draw[->] (fk6.east) to   (gk5.west) ;
  \draw[->] (fk4.east) to   (gk4.west) ;
  \draw[->] (fk5.east) to   (gk5.west) ;
  \draw[->] (fk6.east) to   (gk6.west) ;

 \node[ble2] at (9.5,-1)    (hk1) {$\textbf{8}$};
  \node[ble] at (9.5,-1.6) (hk2) {$6$};
  \node[ble] at (9.5,-2.2) (hk3) {$2$};
  \node[ble] at (9.5,-2.8) (hk4) {$5$};
  \node[ble] at (9.5,-3.4) (hk5) {$1$};
  \node[ble] at (9.5,-4)   (hk6) {$7$};
  \node[ble] at (9.5,-4.6)   (hk7) {$4$};
  \node[ble] at (9.5,-5.2)   (hk8) {$3$};


  \draw[->] (gk1.east) to   (hk2.west) ;
  \draw[->] (gk1.east) to   (hk3.west) ;
  \draw[->] (gk1.east) to   (hk4.west) ;
  \draw[->] (gk1.east) to   (hk5.west) ;
  \draw[->] (gk1.east) to   (hk6.west) ;
  \draw[->] (gk1.east) to   (hk7.west) ;
  \draw[->] (gk1.east) to   (hk8.west) ;
  \draw[->] (gk2.east) to   (hk1.west) ;
  \draw[->] (gk3.east) to   (hk2.west) ;
  \draw[->] (gk4.east) to   (hk3.west) ;
  \draw[->] (gk5.east) to   (hk4.west) ;
  \draw[->] (gk6.east) to   (hk5.west) ;
  \draw[->] (gk7.east) to   (hk6.west) ;
  \draw[->] (gk8.east) to   (hk7.west) ;

  \draw[->] (gk2.east) to   (hk1.west) ;
  \draw[->] (gk3.east) to   (hk2.west) ;
  \draw[->] (gk4.east) to   (hk3.west) ;
  \draw[->] (gk5.east) to   (hk4.west) ;
  \draw[->] (gk6.east) to   (hk5.west) ;
  \draw[->] (gk7.east) to   (hk6.west) ;
  \draw[->] (gk8.east) to   (hk7.west) ;

  \draw[->] (gk2.east) to   (hk2.west) ;
  \draw[->] (gk3.east) to   (hk3.west) ;
  \draw[->] (gk4.east) to   (hk4.west) ;
  \draw[->] (gk5.east) to   (hk5.west) ;
  \draw[->] (gk6.east) to   (hk6.west) ;
  \draw[->] (gk7.east) to   (hk7.west) ;
  \draw[->] (gk8.east) to   (hk8.west) ;
  

 \node[ble] at (11,-1)    (ik1) {$6$};
  \node[ble] at (11,-1.6) (ik2) {$8$};
  \node[ble] at (11,-2.2) (ik3) {$2$};
  \node[ble] at (11,-2.8) (ik4) {$5$};
  \node[ble] at (11,-3.4) (ik5) {$1$};
  \node[ble] at (11,-4)   (ik6) {$7$};
  \node[ble] at (11,-4.6)   (ik7) {$4$};
  \node[ble] at (11,-5.2)   (ik8) {$3$};

  \draw[->] (hk1.east) to   (ik1.west) ;
  \draw[->] (hk1.east) to   (ik2.west) ;
  \draw[->] (hk1.east) to   (ik3.west) ;
  \draw[->] (hk1.east) to   (ik4.west) ;
  \draw[->] (hk1.east) to   (ik5.west) ;
  \draw[->] (hk1.east) to   (ik6.west) ;
  \draw[->] (hk1.east) to   (ik7.west) ;
  \draw[->] (hk1.east) to   (ik8.west) ;
  
  \draw[->] (hk2.east) to   (ik1.west) ;
  \draw[->] (hk3.east) to   (ik2.west) ;
  \draw[->] (hk4.east) to   (ik3.west) ;
  \draw[->] (hk5.east) to   (ik4.west) ;
  \draw[->] (hk6.east) to   (ik5.west) ;
  \draw[->] (hk7.east) to   (ik6.west) ;
  \draw[->] (hk8.east) to   (ik7.west) ;

  \draw[->] (hk2.east) to   (ik2.west) ;
  \draw[->] (hk3.east) to   (ik3.west) ;
  \draw[->] (hk4.east) to   (ik4.west) ;
  \draw[->] (hk5.east) to   (ik5.west) ;
  \draw[->] (hk6.east) to   (ik6.west) ;
  \draw[->] (hk7.east) to   (ik7.west) ;
  \draw[->] (hk8.east) to   (ik8.west) ;


\end{tikzpicture}
\caption{Example run of \textit{top-to-random} shuffle
taken from the key $k_g$ of the length equal to 8 8-value bytes (24-bits).  
Conditioning on the number of steps of the \sst (in this case 8) one can find 
out that: out of the possible $2^{24}$ keys only (exactly) $8!$ keys are 
possible (due to the fact that \sst stopped exactly after $8$ steps) and every 
of $8!$ 
permutations is possible -- moreover each permutation with exactly the 
same probability. \sst leaks bits of the key \ie $k_g[1] = 111$ but does not 
leak any information about the produced permutation.} \label{fig:przebieg-fresh}
\end{center}
\end{figure*}

In  Figure~\ref{fig:przebieg-fresh} a sample   execution of the algorithm is 
presented (for $n=8$). It is an extreme situation in which SST was achieved after $n$ steps 
(the lower bound for the running time) -- please keep in mind that the 
expected running time of this SST is $O(n \log n)$.
 Based on this knowledge 
it is easy to observe that  there is only one possible value of $k_g[1]$ 
provided that \sst lasted $n$ steps -- $k_g[1] = 111 = [8] \setminus [7]$; $k_g[2]$ 
is either $110$ or $111$ so $k_g[2] \in [8] \setminus [6]$ and $k_g[i] \in [n] 
\setminus [n-i]$ for each $i = 1, \ldots, n$. While some  information 
about 
the key $k_g$ leaks from the fact that the \sst stopped just after $n$ steps, 
still no adversary 
can learn anything about the resulting permutation because every permutation is 
generated with exactly the same probability. 

We showed an extreme case 
when SST algorithm for \textsl{Top-To-Random} runs the 
minimal number of steps $n$. Since  $k_g[i] \in [n] \setminus [n-i]$ for $i = 
1, \ldots, n$ so the number of bits of entropy of the secret key $k_g$ adversary 
learns each round $i$ is equal to $\lg i$. 
There are $n!$ possible keys which result in running time $n$ of SST.
An adversary learns a fraction of $n!$ out of $2^{n \lg n} = n^n$ possible 
keys. Moreover it can predict a prefix of the $k_g$ with high probability. 
\subsubsection{Riffle Shuffle}

Consider the algorithm \ksa{${}^{**}_\texttt{RS,ST}$} with the shuffling 
procedure 
corresponding to (time reversal of) \textsl{Riffle Shuffle}. Recall the 
shuffling:
at each step  assign each card a random bit, then put all cards with bit 
assigned 0 to the top keeping the relative ordering.
Sample execution is presented in Fig. \ref{fig:riffle_przebieg}.

\tikzstyle{RectObject}=[rectangle,fill=white,draw,line width=0.5mm]
\tikzstyle{line}=[draw]
\tikzstyle{arrow}=[draw, -latex]

\tikzstyle{block} = [draw,fill=gray!40,minimum size=0.5em]
\tikzstyle{ble} = [draw,minimum size=0.5em]
\tikzstyle{ble2} = [draw,fill=red!20,minimum size=0.5em]

\tikzstyle{block2} = [draw,fill=gray!90,minimum size=0.5em]
\def\radius{.7mm} 
\tikzstyle{branch}=[fill,shape=circle,minimum size=3pt,inner sep=0pt]

\begin{figure*}
\begin{center}
\begin{tikzpicture}[>=latex',scale=.9,every node/.style={scale=.9}]

   \node  at (0.75,-0.2) {step: 1};
   
   \node  at (3.75,-0.2) {step: 2};
   \node  at (6.75,-0.2) {step: 3};

  \node[block] at (-1,-1)   (ak1) {$1$};
  \node[block] at (-1,-1.6) (ak2) {$2$};
  \node[block] at (-1,-2.2) (ak3) {$3$};
  \node[block] at (-1,-2.8) (ak4) {$4$};
  \node[block] at (-1,-3.4) (ak5) {$5$};
  \node[block] at (-1,-4)   (ak6) {$6$};
  \node[block] at (-1,-4.6)   (ak7) {$7$};
  \node[block] at (-1,-5.2)   (ak8) {$8$};

  \node[ble] at (-0.5,-1)   (ak1b) {$0$};
  \node[ble] at (-0.5,-1.6) (ak2b) {$0$};
  \node[ble] at (-0.5,-2.2) (ak3b) {$0$};
  \node[ble] at (-0.5,-2.8) (ak4b) {$0$};
  \node[ble] at (-0.5,-3.4) (ak5b) {$1$};
  \node[ble] at (-0.5,-4)   (ak6b) {$1$};
  \node[ble] at (-0.5,-4.6)   (ak7b) {$1$};
  \node[ble] at (-0.5,-5.2)   (ak8b) {$1$};

 \node[block] at (2,-1)    (bk1) {$1$};
  \node[block] at (2,-1.6) (bk2) {$2$};
  \node[block] at (2,-2.2) (bk3) {$3$};
  \node[block] at (2,-2.8) (bk4) {$4$};
  \node[block] at (2,-3.4) (bk5) {$5$};
  \node[block] at (2,-4)   (bk6) {$6$};
  \node[block] at (2,-4.6)   (bk7) {$7$};
  \node[block] at (2,-5.2)   (bk8) {$8$};

  \node[ble] at (2.5,-1)   (bk1b) {$0$};
  \node[ble] at (2.5,-1.6) (bk2b) {$0$};
  \node[ble] at (2.5,-2.2) (bk3b) {$1$};
  \node[ble] at (2.5,-2.8) (bk4b) {$1$};
  \node[ble] at (2.5,-3.4) (bk5b) {$0$};
  \node[ble] at (2.5,-4)   (bk6b) {$0$};
  \node[ble] at (2.5,-4.6)   (bk7b) {$1$};
  \node[ble] at (2.5,-5.2)   (bk8b) {$1$};

  \draw[->] (ak1b.east) to   (bk1.west) ;
  \draw[->] (ak2b.east) to   (bk2.west) ;
  \draw[->] (ak3b.east) to   (bk3.west) ;
  \draw[->] (ak4b.east) to   (bk4.west) ;
  \draw[->] (ak5b.east) to   (bk5.west) ;
  \draw[->] (ak6b.east) to   (bk6.west) ;
  \draw[->] (ak7b.east) to   (bk7.west) ;
  \draw[->] (ak8b.east) to   (bk8.west) ;

 \node[block] at (5,-1)   (ck1) {$1$};
  \node[block] at (5,-1.6) (ck2) {$2$};
  \node[block] at (5,-2.2) (ck3) {$5$};
  \node[block] at (5,-2.8) (ck4) {$6$};
  \node[block] at (5,-3.4) (ck5) {$3$};
  \node[block] at (5,-4)   (ck6) {$4$};
  \node[block] at (5,-4.6)   (ck7) {$7$};
  \node[block] at (5,-5.2)   (ck8) {$8$};

  \node[ble] at (5.5,-1)   (ck1b) {$0$};
  \node[ble] at (5.5,-1.6) (ck2b) {$1$};
  \node[ble] at (5.5,-2.2) (ck3b) {$0$};
  \node[ble] at (5.5,-2.8) (ck4b) {$1$};
  \node[ble] at (5.5,-3.4) (ck5b) {$0$};
  \node[ble] at (5.5,-4)   (ck6b) {$1$};
  \node[ble] at (5.5,-4.6)   (ck7b) {$0$};
  \node[ble] at (5.5,-5.2)   (ck8b) {$1$};

    \draw[->] (bk1b.east) to   (ck1.west) ;
  \draw[->] (bk2b.east) to   (ck2.west) ;
  \draw[->] (bk3b.east) to   (ck5.west) ;
  \draw[->] (bk4b.east) to   (ck6.west) ;
  \draw[->] (bk5b.east) to   (ck3.west) ;
  \draw[->] (bk6b.east) to   (ck4.west) ;
  \draw[->] (bk7b.east) to   (ck7.west) ;
  \draw[->] (bk8b.east) to   (ck8.west) ;

 \node[block] at (8,-1)   (dk1) {$1$};
  \node[block] at (8,-1.6) (dk2) {$5$};
  \node[block] at (8,-2.2) (dk3) {$3$};
  \node[block] at (8,-2.8) (dk4) {$7$};
  \node[block] at (8,-3.4) (dk5) {$2$};
  \node[block] at (8,-4)   (dk6) {$6$};
  \node[block] at (8,-4.6)   (dk7) {$4$};
  \node[block] at (8,-5.2)   (dk8) {$8$};

  \draw[->] (ck1b.east) to   (dk1.west) ;
  \draw[->] (ck3b.east) to   (dk2.west) ;
  \draw[->] (ck5b.east) to   (dk3.west) ;
  \draw[->] (ck7b.east) to   (dk4.west) ;
  
  \draw[->] (ck2b.east) to   (dk5.west) ;
  \draw[->] (ck4b.east) to   (dk6.west) ;
  \draw[->] (ck6b.east) to   (dk7.west) ;
  \draw[->] (ck8b.east) to   (dk8.west) ;

\end{tikzpicture}
\caption{Example run of  (time reversal of) \textsl{Riffle Shuffle} of $n=8$ 
cards and key $k_g=[00001111,00110011,01010101,\ldots]$.
At each step the left column represents current permutation, whereas the right 
one currently assigned bits. 
See description in the text.
 } \label{fig:riffle_przebieg}
\end{center}
\end{figure*}
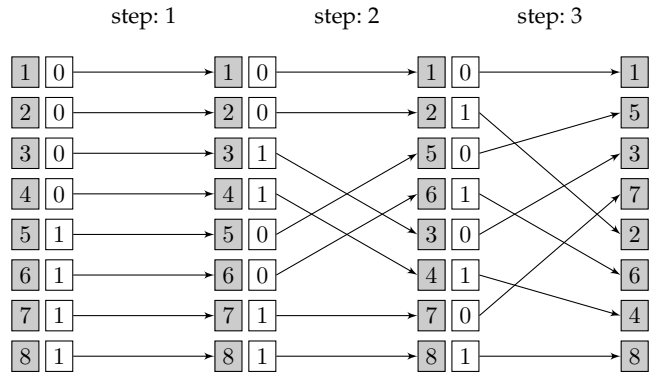 
\noindent Recall the construction of SST for this example:
\noindent 
\begin{quotation}
\noindent 
Consider, initially \textsl{unmarked},  all ${n\choose 2}$ pairs of cards; at 
each step \textsl{mark} pair $(i,j)$ if cards $i$ and 
$j$ were assigned different bits.  
\end{quotation}
Assume for simplicity that $n$ is a power of 2. 
How many \textsl{pairs} can be marked in first step? For example, if there is 
only one bit equal to 0 and all other 
are equal to 1, then $n-1$ pairs are marked. It is easy to see that the maximal number of
pairs will be marked if exactly half (\ie ${n\over 2}$)
cards were assign bit 1 and other half were assign bit 0. Then ${n^2\over 4}$ 
pairs will be marked (all pairs of form 
$(i,j)$: $i<j$ and cards $i$ and $j$ were assigned different bits). After such 
a step there will be two groups of cards: the ones assigned 0 and ones 
assigned 1. All pairs of cards within each group are not marked. Then again, to 
mark as many pairs as possible exactly half of the cards 
in each group must be assigned bit 0 and other half bit 1 -- then ${n^2\over 8}$ 
new pairs will be marked. Iterating this procedure,
we see that the shortest execution takes $\lg n$ steps after which all 
$${n^2\over 4}+{n^2\over 8}+{n^2\over 16} + \ldots +  {n^2\over 2^{1+\lg n}} 
={n\choose 2}$$
pairs are marked. In Figure \ref{fig:riffle_przebieg} such sample execution is 
presented.
There are exactly $n!$ keys resulting in $\lg n$ steps of SST (each resulting in 
different permutation). Thus 
the adversary learns a tiny fraction of possible keys: $n!$ out of $2^{n\lg 
n} = n^n$.

\subsubsection{Summary}
As we have seen the choice of shuffle algorithm matters for  \ksa{${}^{**}$}. In 
any case the resulting permutation is not revealed at all,
but some knowledge on key can be obtained. In both presented examples  an adversary who 
observed that the time to generate a permutation was minimal (for a given 
SST) learns that one out of $n!$ keys was used (out of $n^n$ possible keys). 
However there is a huge difference between these two shuffling schemes:
\begin{itemize}
 \item in \textsl{Riffle Shuffle} an adversary learns that in the first 
step one key out of ${n \choose n/2}$ possible keys (leading to the shortest 
running time) was used. For the second step there are exactly ${n/2 \choose 
n/4}^2$ for each choice of the key made in the fist step, and similarly for the 
next steps. This makes recovering any information about the key used infeasible.
\item in \textsl{Top-To-Random} the situation of an adversary is quite 
different. It gets all $\lg n$ bits of the key used in the 
first step of the algorithm (provided the algorithm run for $n$ steps). In the 
$i$th round it gets information about $\lg n-i+1$ bits.
\end{itemize}
Recall that these were the extreme and very unlikely cases. In both cases the 
probability of the shortest run is equal to $n!/n^n$ ($=7.29 \times 10^{-55}$ 
for $n = 128$).

Another advantage of \textsl{Riffle Shuffle} is that on average the number of 
required bits is equal to $2 n \lg n$ which is only $2$ times more than the 
minimal number of bits (as presented in the example). For 
\textsl{Top-To-Random} on average the number of bits required is equal to $n 
\log n \lg n$ which is $\log n$ times larger than the minimal number of bits.
This follows from the fact that the variance of SST for \textsl{Riffle Shuffle} 
is much smaller than the variance for \textsl{Top-To-Random}. Smaller variance 
means that fewer information can be gained from the running time.
While the expected number of bits required by \ksa{${}^{**}$}  may seem to be 
high, for the post-quantum algorithms like McBits it is not an issue  -- McBits 
requires keys of size $2^{16} - 2^{20}$ bytes, \textsl{Top-To-Random} needs 
around 12,000 bits and \textsl{Riffle Shuffle} needs on average $4096$ bits to 
generate a uniform permutation for $n = 256$.

It is also worth mentioning that in real-life systems we could 
artificially additionally mask the duration of the algorithm:
 If we know values of $ET$ for given SST we can run it always for at least $ET$ steps (similarly, to avoid extremely 
 long execution we can let at most $ET+c\cdot VarT$ steps provided $VarT$ is known). 

\section{Conclusions}\label{sec:conclusions}
We presented the benefits of using Strong Stationary Times  in cryptographic 
schemes (pseudo random permutation generators). These algorithms   have a 
``health-check'' built-in and 
guarantee the best possible properties (when it comes to the quality of 
randomness of the resulting permutation). 

We showed how SST algorithms may be used to construct timing attack immune 
implementations out of timing attack susceptible ones. 
The SST-based masking leads to comparable (to constant-time 
implementation) slow-down (but it is platform independent), see 
Figure~\ref{tab:efficiency} for efficiency comparison.

We showed that algorithms using SST achieve 
better 
security guarantees than any algorithm which runs predefined number of steps.

We also proved a better bound for the mixing-time of the \ctrt 
shuffling process 
 which is 
used in \rcc and showed that different, more efficient 
shuffling methods (\ie time reversal of \riffle) may be used as \ksa. This 
last observation shows that the gap between mathematical model (4096 bits 
required) and reality (2048 allowed as maximum length of \rcc) is not that big 
as previously thought (bound of 23~037 by Mironov~\cite{Mironov2002}).


  \bibliographystyle{plain}
 \bibliography{library}

\begin{thebibliography}{10}

\bibitem{Intel_perm}
{Intel 64 and IA-32 Architectures. Software Developer's Manual}, 2016.

\bibitem{albrecht2015lucky}
Martin~R. Albrecht and Kenneth~G. Paterson.
\newblock {Lucky Microseconds: A Timing Attack on Amazon's s2n Implementation
  of TLS}.
\newblock {\em Advances in Cryptology - EUROCRYPT}, 9665:622----643, 2016.

\bibitem{Aldous1986}
David Aldous and Persi Diaconis.
\newblock {Shuffling cards and stopping times}.
\newblock {\em American Mathematical Monthly}, 93(5):333--348, 1986.

\bibitem{Aldous1987}
David Aldous and Persi Diaconis.
\newblock {Strong Uniform Times and Finite Random Walks}.
\newblock {\em Advances in Applied Mathematics}, 97:69--97, 1987.

\bibitem{rc4royal}
Nadhem AlFardan, Daniel~J Bernstein, Kenneth~G Paterson, Bertram Poettering,
  and Jacob C~N Schuldt.
\newblock {On the Security of RC4 in TLS}.
\newblock In {\em Presented as part of the 22nd USENIX Security Symposium
  (USENIX Security 13)}, pages 305--320, Washington, D.C., 2013. USENIX.

\bibitem{BenesNetwork}
V.E. Bene\v{s}.
\newblock {\em {Mathematical theory of connecting networks and telephone
  traffic}}.
\newblock Academic Press, New York :, 1965.

\bibitem{McBits}
Daniel~J Bernstein, Tung Chou, and Peter Schwabe.
\newblock {McBits: fast constant-time code-based cryptography}.
\newblock In {\em Cryptographic Hardware and Embedded Systems - CHES 2013},
  2013.

\bibitem{Brumley2005}
David Brumley and Dan Boneh.
\newblock {Remote timing attacks are practical}.
\newblock {\em Computer Networks}, 48(5):701--716, 2005.

\bibitem{Chaum1983}
David Chaum.
\newblock {Blind Signatures for Untraceable Payments}.
\newblock In {\em Advances in Cryptology}, pages 199--203. Springer US, Boston,
  MA, 1983.

\bibitem{Dhem1998}
Jean-Fran\c{c}ois Dhem, Fran\c{c}ois Koeune, Philippe-Alexandre Leroux, Patrick
  Mestr\'{e}, Jean-Jacques Quisquater, and Jean-Louis Willems.
\newblock {A Practical Implementation of the Timing Attack}.
\newblock {\em Smart Card. Research and Applications Third International
  Conference, CARDIS98}, 1820:167--182, 1998.

\bibitem{Diaconis1981a}
Persi Diaconis and Mehrdad Shahshahani.
\newblock {Generating a random permutation with random transpositions}.
\newblock {\em Zeitschrift fur Wahrscheinlichkeitstheorie und Verwandte
  Gebiete}, 57(2):159--179, 1981.

\bibitem{Fill1998}
James~Allen Fill.
\newblock {An interruptible algorithm for perfect sampling via Markov chains}.
\newblock {\em The Annals of Applied Probability}, 8(1):131--162, February
  1998.

\bibitem{Fluhrer2001a}
S~Fluhrer, I~Mantin, and A~Shamir.
\newblock {Weaknesses in the key scheduling algorithm of RC4}.
\newblock {\em Selected areas in cryptography}, 2001.

\bibitem{Fluhrer2000}
Scott~R Fluhrer and David~a. McGrew.
\newblock {Statistical Analysis of the Alleged RC4 Keystream Generator}.
\newblock {\em Fast Software Encryption, 7th International Workshop}, pages
  19--30, 2000.

\bibitem{Golic1997}
J~Golic.
\newblock {Linear Statistical Weakness of Alleged RC4 Keystream Generator}.
\newblock In Walter Fumy, editor, {\em Advances in Cryptology — EUROCRYPT
  ’97}, volume 1233 of {\em Lecture Notes in Computer Science}. Springer
  Berlin Heidelberg, Berlin, Heidelberg, July 1997.

\bibitem{Kocher1996}
Paul~C. Kocher.
\newblock {Timing Attacks on Implementations of Diffie-Hellman, RSA, DSS, and
  Other Systems}.
\newblock pages 104--113. Springer, Berlin, Heidelberg, 1996.

\bibitem{Zagorski_RST}
Michal Kulis, Paweł Lorek, and Filip Zag\'{o}rski.
\newblock {Randomized stopping times and provably secure pseudorandom
  permutation generators}.
\newblock In {\em Phan RW., Yung M. (eds) Paradigms in Cryptology – Mycrypt
  2016. Malicious and Exploratory Cryptology. Mycrypt 2016. Lecture Notes in
  Computer Science, vol 10311. Springer}, pages 145----167, 2017.

\bibitem{LorekMarkowski2014}
Paweł Lorek and Piotr Markowski.
\newblock {Monotonicity requirements for efficient exact sampling with Markov
  chains.}
\newblock {\em To appear in Markov Processes And Related Fields}, 2017.

\bibitem{ManSha2001}
Itsik Mantin and Adi Shamir.
\newblock {A Practical Attack on Broadcast RC4}.
\newblock {\em Fast Software Encryption, 8th International Workshop, Yokohama,
  Japan}, pages 152--164, 2001.

\bibitem{Matthews1988}
Peter Matthews.
\newblock {A strong uniform time for random transpositions}.
\newblock {\em J. Theoret. Probab.}, 1(4):411--423, 1988.

\bibitem{Mironov2002}
Ilya Mironov.
\newblock {(Not So) Random Shuffles of RC4}.
\newblock {\em Advances in Cryptology—CRYPTO 2002}, 2002.

\bibitem{Mossel2004}
Elchanan Mossel, Yuval Peres, and Alistair Sinclair.
\newblock {Shuffling by semi-random transpositions}.
\newblock {\em Foundations of Computer Science}, pages 572--581, 2004.

\bibitem{Pereida}
Cesar Pereida, Garc\'{\i}a Billy, and Bob Brumley.
\newblock {Constant-Time Callees with Variable-Time Callers}.
\newblock Technical report.

\bibitem{PereidaGarcia2016}
Cesar {Pereida Garc\'{\i}a}, Billy~Bob Brumley, and Yuval Yarom.
\newblock {Make Sure DSA Signing Exponentiations Really are Constant-Time}.
\newblock pages 1639--1650, 2016.

\bibitem{Pornin}
Thomas Pornin.
\newblock {BearSSL - Constant-Time Crypto}.

\bibitem{PorninThomas}
{Pornin Thomas}.
\newblock {BearSSL - Constant-Time Mul}.

\bibitem{Propp1996}
James~Gary Propp and David~Bruce Wilson.
\newblock {Exact sampling with coupled Markov chains and applications to
  statistical mechanics}.
\newblock {\em Random Structures and Algorithms}, 9:223--252, 1996.

\bibitem{Spritz}
Jacob C. N. Schuldt; Ronald~L. Rivest.
\newblock {Spritz—a spongy RC4-like stream cipher and hash function}.
\newblock Technical report, 2014.

\bibitem{rc4a}
Bart~Preneel {Souradyuti Paul}.
\newblock {A New Weakness in the RC4 Keystream Generator and an Approach to
  Improve the Security of the Cipher}.
\newblock In Bimal Roy and Willi Meier, editors, {\em Fast Software
  Encryption}, volume 3017 of {\em Lecture Notes in Computer Science}. Springer
  Berlin Heidelberg, Berlin, Heidelberg, 2004.

\bibitem{rc4plus}
Bart~Preneel {Souradyuti Paul}.
\newblock {A New Weakness in the RC4 Keystream Generator and an Approach to
  Improve the Security of the Cipher}.
\newblock In Bimal Roy and Willi Meier, editors, {\em Fast Software
  Encryption}, volume 3017 of {\em Lecture Notes in Computer Science}. Springer
  Berlin Heidelberg, Berlin, Heidelberg, 2004.

\bibitem{Cachebleed}
Yuval Yarom, Daniel Genkin, and Nadia Heninger.
\newblock {CacheBleed : A Timing Attack on OpenSSL Constant Time RSA}.
\newblock {\em CHES}, 2016.

\bibitem{vmpc}
Bartosz Zoltak.
\newblock {VMPC One-Way Function and Stream Cipher}.
\newblock In {\em Fast Software Encryption}, pages 210--225. 2004.

\end{thebibliography}

 \end{document}